\documentclass[12pt,a4paper]{article}
\usepackage{amsmath}
\usepackage{amsfonts}
\usepackage{amssymb}
\usepackage[utf8]{inputenc}
\usepackage[T1]{fontenc}
\usepackage{amssymb,amsmath,amsthm}
\usepackage{graphicx}
\usepackage{enumerate}
\usepackage{longtable}
\date{}

\DeclareMathOperator{\conv}{conv}
\DeclareMathOperator{\bd}{bd}
\DeclareMathOperator{\lin}{lin}

\DeclareMathOperator{\sgn}{sgn}

\DeclareMathOperator{\inte}{int}
\DeclareMathOperator{\as}{as}
\newtheorem{twr}{Theorem}[section]
\newtheorem{lem}[twr]{Lemma}

\theoremstyle{remark}
\newtheorem{remark}[twr]{Remark}

\begin{document}

\title{On the Banach-Mazur Distance in Small Dimensions}

\author{Tomasz Kobos\\
Faculty of Mathematics and Computer Science\\
Jagiellonian University in Cracow\\
\small Lojasiewicza 6, 30-348 Krakow, Poland\\
\small\texttt{tomasz.kobos@uj.edu.pl}\\
\and 
Marin Varivoda\\
Faculty of Science - Department of Mathematics\\ 
University of Zagreb\\ 
\small Bijenička cesta 30, 10000 Zagreb, Croatia\\
\small {\texttt{mvarivod.math@pmf.hr}}
}

\maketitle

\begin{abstract}
We establish some results on the Banach-Mazur distance in small dimensions. Specifically, we determine the Banach-Mazur distance between the cube and its dual (the cross-polytope) in $\mathbb{R}^3$ and $\mathbb{R}^4$. In dimension three this distance is equal to $\frac{9}{5}$, and in dimension four, it is equal to 2. These findings confirm well-known conjectures, which were based on numerical data. Additionally, in dimension two, we use the asymmetry constant to provide a geometric construction of a family of convex bodies that are equidistant to all symmetric convex bodies.
\end{abstract}

\footnotetext{
\textbf{Keywords} Banach-Mazur Distance, Hypercube, Cross-Polytope, Asymmetry Constant\\
\textbf{MSC Classification} 52A40, 52A10, 52A27}     

\section{Introduction}

Let $n \geq 2$ be an integer. A \textit{convex body} in $\mathbb{R}^n$ is a compact, convex set with a non-empty interior. A convex body will be called \textit{centrally-symmetric} (or just \textit{symmetric}) if it has a center of symmetry. For two convex bodies $K, L \subseteq \mathbb{R}^n$ we define their \textit{Banach-Mazur} distance as 
$$d_{BM}(K,L) = \inf \{r >0 : K+ u \subseteq T(L+v) \subseteq r(K+u) \},$$ 
where the infimum is taken over all invertible linear operators $T : \mathbb{R}^n \to \mathbb{R}^n$ and vectors $u, v \in \mathbb{R}^n$. One can easily check that this infimum is attained by some operator. Moreover, if $K$ and $L$ are both symmetric with respect to the origin, then it is attained for $u = v = 0$. We note that this is a multiplicative distance when it is considered as a distance between equivalence classes of convex bodies -- the distance between a convex body and its non-degenerate affine copy is by definition equal to $1$. 

The Banach-Mazur distance is a well-established notion of functional analysis, as the Banach-Mazur distance between unit balls of two norms in $\mathbb{R}^n$ can be interpreted as the distance between two $n$-dimensional normed spaces. This was actually the original definition of this notion that was introduced by Banach in \cite{banach}. One can say that the Banach-Mazur distance serves the purpose of comparing the geometric properties of two normed spaces and quantifies how essentially different the spaces are. This is reflected in its numerous important applications in the fields of convex geometry, discrete geometry and local theory of Banach spaces. This notion has already been extensively studied by many authors, leading to some remarkable results. One very famous example is the Gluskin construction \cite{gluskin} of symmetric random polytopes in $\mathbb{R}^n$ with the Banach-Mazur distance of order $cn$. Random construction of Gluskin was a major breakthrough in the local theory of Banach spaces, as the method turned out to have many more possible applications and consequently had a profound impact on this field. An excellent reference is the monograph of N. Tomczak-Jaegermann \cite{tomczak} which in large part is devoted to a detailed study of the Banach-Mazur distance from the viewpoint of functional analysis.

It should be emphasized however, that the vast majority of established results regarding the Banach-Mazur distance are asymptotic in nature. In other words, these results mostly describe the behavior of the Banach-Mazur distance as the dimension tends to infinity. On the other hand, the non-asymptotic properties of the Banach-Mazur distance seem to be quite elusive, and even in very small dimensions they are surprisingly difficult to establish. For example, it is known that the maximal possible distance between two symmetric bodies in $\mathbb{R}^n$ is asymptotically of order $cn$ (which follows from John's Ellipsoid Theorem and Gluskin's random construction of convex bodies), but the precise value of this maximal distance is known only for $n=2$. In this case Stromquist \cite{stromquist} proved that the distance between the square and the regular hexagon is equal to $\frac{3}{2}$, and this is the maximal possible distance between a pair of planar symmetric convex bodies. Actually, there are rather few situations in which the Banach-Mazur distance between a pair of convex bodies has been determined precisely. One example illustrating this difficulty is the case of the cube and the cross-polytope (regular octahedron) in $\mathbb{R}^3$. These are perhaps the two simplest symmetric convex polytopes in the three-dimensional space, and yet their Banach-Mazur distance was not determined. Numerical results suggested that this distance is equal to $\frac{9}{5}$, and that the corresponding distance in dimension four is equal to $2$ (for a reference about the numerical data see paper \cite{xue} by Xue). In the planar case this distance is obviously equal to $1$, and in the asymptotic setting it has been known for a long time that $d_{BM}(\mathcal{C}_n, \mathcal{C}^*_n)$ is of the order $\sqrt{n}$, where by $\mathcal{C}_n$ and $\mathcal{C}^*_n$ we denote the unit cube and its dual (the cross-polytope) in $\mathbb{R}^n$ respectively. More precisely, there exist absolute constants $c_1, c_2 > 0$ such that 
$$c_1 \sqrt{n} \leq d_{BM}(\mathcal{C}_n, \mathcal{C}^*_n) \leq c_2 \sqrt{n}$$
for every $n \geq 1$ (see for example Proposition 37.6 in \cite{tomczak}). However, these asymptotic estimates do not say a lot about the small dimensional cases. Some specific upper and lower bounds for the Banach-Mazur distance between the $n$-dimensional cube and the cross-polytope were given by Xue in \cite{xue}.

 It is worth noting that even more generally, the maximal possible distance of a symmetric convex body to the $n$-dimensional cube (or the cross-polytope) has been studied by several authors, but also mainly with a focus on asymptotic properties (see for example: \cite{bourgain}, \cite{giannopoulos}, \cite{tikhomirov}, \cite{youssef}). In small dimensions, the best possible upper bound for the maximal possible distance of a symmetric convex body to the cube was given by Taschuk in \cite{taschuk}. However, for determining the distance between the three and four dimensional cube and the cross-polytope, the main difficulty lies in establishing the lower bound. In this case, it is not difficult to find linear operators that provide the upper bound of $\frac{9}{5}$ and $2$, respectively. 

 The main goal of this paper is to establish some results in small dimensions, in which the Banach-Mazur distance can be determined precisely. The paper is divided into three distinct sections, each dealing with a different dimension. Section \ref{sect3d} is concerned with the three-dimensional case, where we give a geometric proof of the fact that the Banach-Mazur distance between the cube and the cross-polytope is equal to $\frac{9}{5}$, hence confirming the well-known conjecture (Theorem \ref{twr3d}). Our approach is based on a simple two-dimensional lemma, which somewhat explains the role played by the number $\frac{9}{5}$ (see Lemma \ref{lemkwadrat}). Moreover, we are able to characterize all linear operators that achieve equality.

In Section \ref{sect4d}, we consider the same question in dimension four. In this case, we prove that the Banach-Mazur distance between the four-dimensional cube and its dual (the cross-polytope or the unit ball of the $\ell_1$ norm in $\mathbb{R}^4$) is equal to $2$, again confirming the result suggested by the numerical data. However, our approach is completely different than in the three-dimensional case and involves a detailed combinatorial analysis. In Remark \ref{remark} we provide an additional observation related to the $n$-dimensional case and the best constant $c$ in the inequality $d_{BM}(\mathcal{C}_n, \mathcal{C}^*_n) \geq c \sqrt{n}$ that is currently known.

In Section \ref{sect2d}, we move on to dimension two and give a geometric construction of a family of planar convex bodies with some special metric properties. The $n$-dimensional simplex is a convex body well-known for its numerous remarkable features. It has been extensively studied, also from the point of view of Banach-Mazur distance. (see for example: \cite{fleicher}, \cite{gordon}, \cite{jimenez}, \cite{kobos}, \cite{lassak}, \cite{novotny}, \cite{palmon}). We shall focus on its following well-known and interesting property: it is equidistant to all symmetric convex bodies, with the distance being equal to $n$ (see for example \cite{grunbaum} or Corollary 5.8 in \cite{gordon}). Moreover, it is known that the simplex is the unique convex body with this property. It is therefore natural to ask if the simplex is the unique convex body that is equidistant to all symmetric convex bodies (not necessarily with the distance equal to $n$)? It turns out that in the planar case the answer is negative. For all $r \in (\frac{7}{4}, 2)$ we prove the existence of continuum many affinely non-equivalent convex pentagons $K \subseteq \mathbb{R}^2$ satisfying $d(K,L)=r$ for every symmetric convex body $L \subseteq \mathbb{R}^2$ (Theorem \ref{twr2d}). To do this, we rely heavily on the properties of a classical affine invariant of convex bodies -- the asymmetry constant. This common distance $r$ is exactly the asymmetry constant of $K$. It is worth emphasizing, that by using the asymmetry constant, we are able to determine the Banach-Mazur distance between a large number of pairs of convex bodies in one go. This is a rather unusual situation, as each of the two preceding sections is devoted to determining the distance only for a specific pair. We note that all methods employed in the paper can be considered to be completely elementary. 

Throughout the paper by $\| \cdot \|_{\infty}$ we will denote the maximum norm in $\mathbb{R}^n$.

\section{Banach-Mazur distance between the cube and the cross-polytope in the three-dimensional case}
\label{sect3d}

\begin{figure}
    \centering
    \includegraphics[scale=0.15]{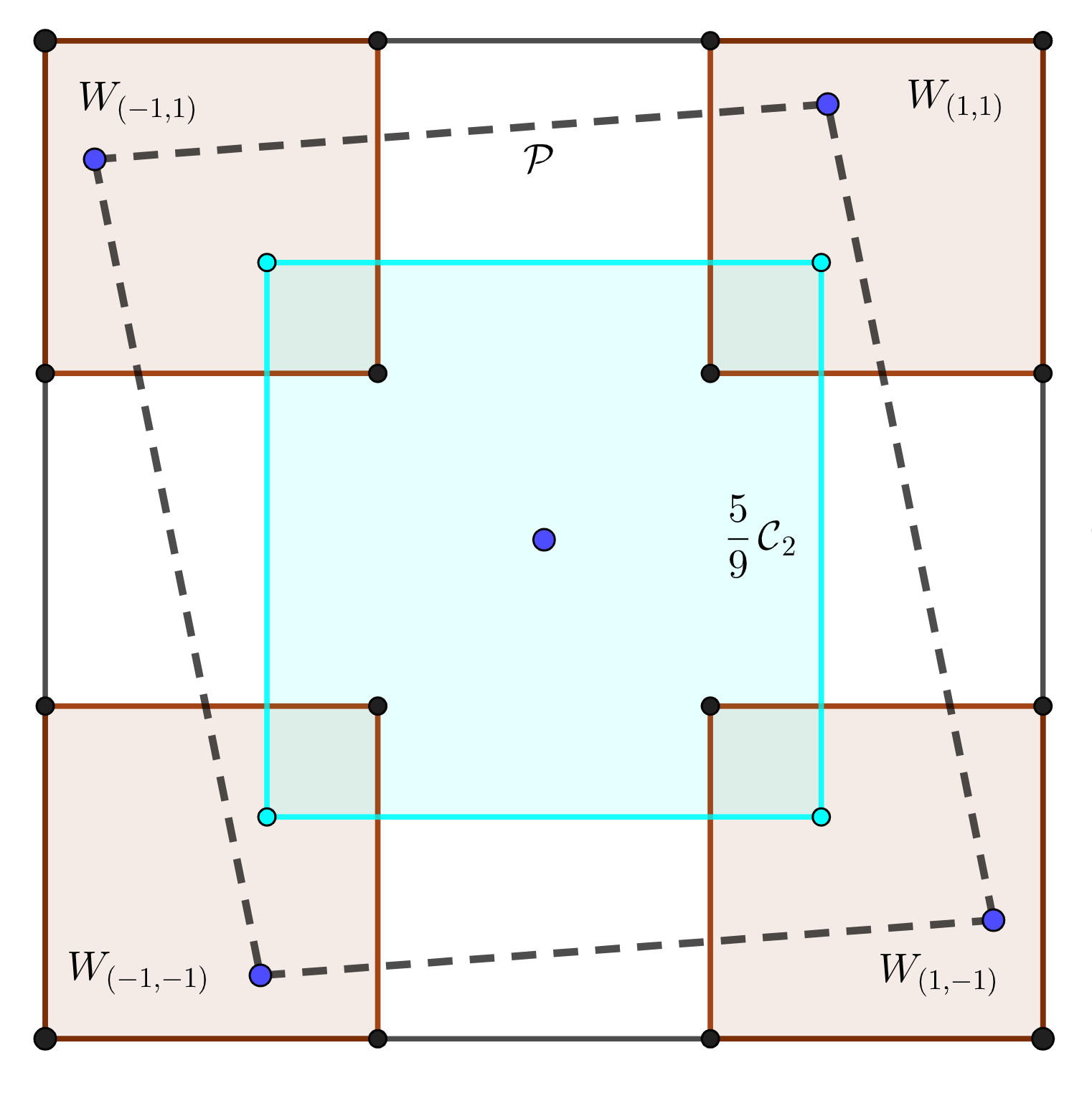}
    \caption{Illustration of Lemma \ref{lemkwadrat}. A $0$-symmetric parallelogram $\mathcal{P}$ is contained in the square $\mathcal{C}_2$ but contains a smaller square $\frac{5}{9} \mathcal{C}_2$. In this case, each of the $4$ squares $W_{(\varepsilon_1, \varepsilon_2)}$ has to contain exactly one vertex of $\mathcal{P}.$}
    \label{ryskwadrat1}
\end{figure}

In this section, we determine the distance between the three-dimensional cube $\mathcal{C}_3$ and the cross-polytope (regular octahedron) $\mathcal{C}_3^*$, providing the positive answer to the conjecture of $\frac{9}{5}$. In order to do this, we will use orthogonal projections onto certain two-dimensional subspaces. The following simple two-dimensional lemma, that is established by means of elementary geometry, represents the reduction of the three-dimensional case to the two-dimensional problem. It can be easily seen that $\frac{5}{9}$ can not be replaced by a smaller number here. See Figure \ref{ryskwadrat1} for illustration. 

\begin{lem}
\label{lemkwadrat}
Let $\mathcal{P} \subseteq \mathbb{R}^2$ be a $0$-symmetric parallelogram in the plane such that $\frac{5}{9}\mathcal{C}_2 \subseteq \mathcal{P} \subseteq \mathcal{C}_2$. For $\varepsilon_1, \varepsilon_2 \in \{-1, 1\}$ let $W_{(\varepsilon_1, \varepsilon_2)} \subseteq \mathbb{R}^2$ be a square defined as
$$W_{(\varepsilon_1, \varepsilon_2)}= \left \{ (x,y) \in \mathbb{R}^2 : \frac{1}{3} \leq \varepsilon_1x, \varepsilon_2y, \leq 1 \right \}.$$
Then each of the $4$ squares $W_{(\varepsilon_1, \varepsilon_2)}$ (where $\varepsilon_1, \varepsilon_2 \in \{-1, 1\}$) contains exactly one vertex of $\mathcal{P}$.
\end{lem}

\begin{figure}
    \centering
    \includegraphics[scale=0.15]{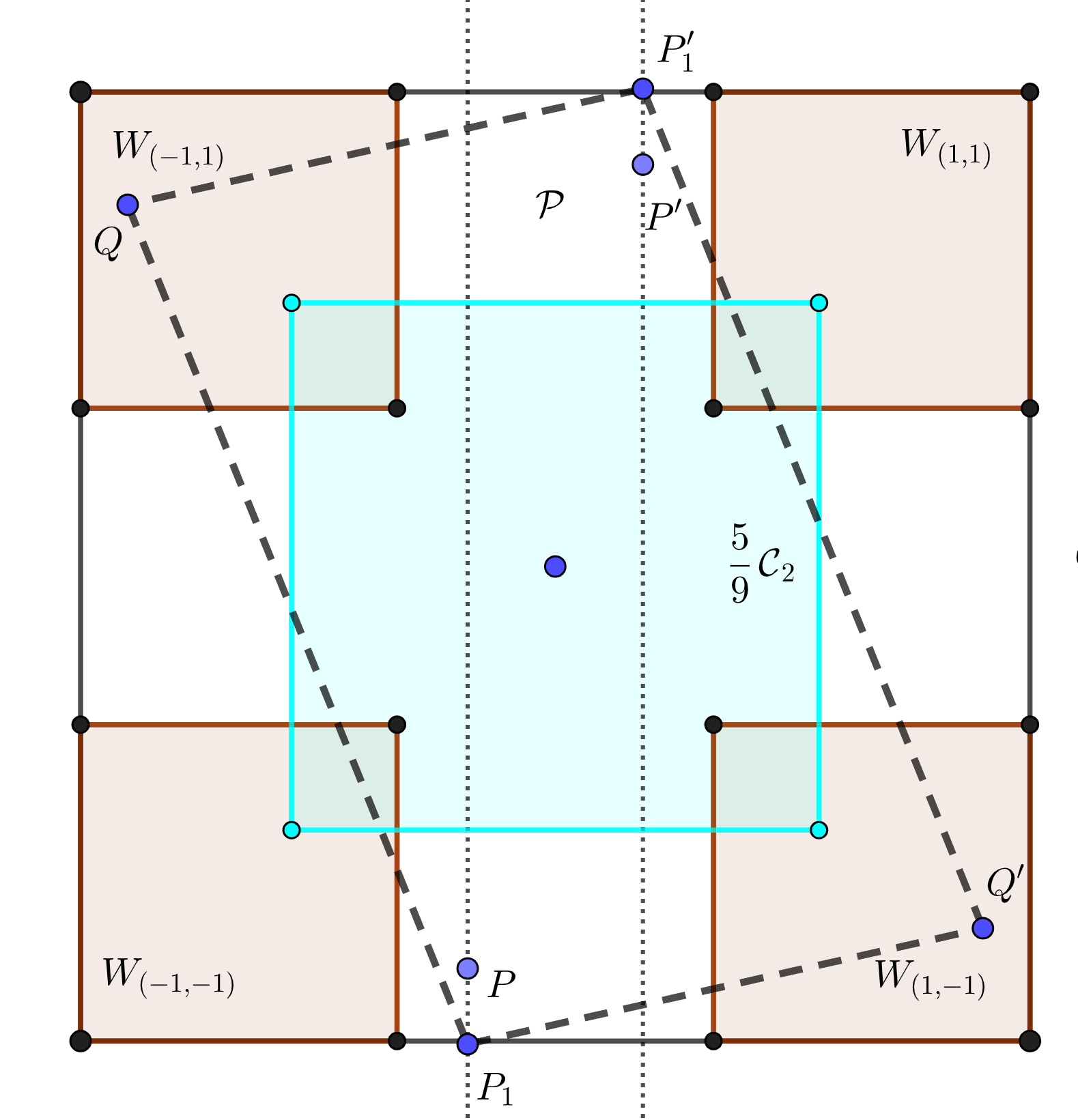}
    \caption{Proof of Lemma \ref{lemkwadrat}. If we assume that a parallelogram $\mathcal{P}$ does not satisfy the desired condition, then by projecting $P$ onto a side of $\mathcal{C}_2$ we get a larger $0$-symmetric parallelogram, that is contained in $\mathcal{C}_2$ but still does not satisfy the desired condition. Hence, we can suppose that $P$ is on a side of $\mathcal{C}_2$.}
    \label{ryskwadrat2}
\end{figure}

\emph{Proof.} We start by proving that each vertex of $\mathcal{P}$ belongs to some square $W_{(\varepsilon_1, \varepsilon_2)}$. Let us assume that $P, Q, P', Q'$ are vertices of $\mathcal{P}$ and some vertex $P$ of $\mathcal{P}$ does not belong to any of the squares $W_{(\varepsilon_1, \varepsilon_2)}$. We can suppose that the situation is like in the Figure \ref{ryskwadrat2}.

If $P_1$ is the projection of $P$ to the corresponding side of $\mathcal{C}_2$, then the segment $QQ'$ cuts the line $PP_1$. Therefore, the parallelogram $P_1QP_1'Q'$ contains $\mathcal{P}$. Thus, we can assume that $P$ is on the boundary of $\mathcal{C}_2$.

\begin{figure}
    \centering
    \includegraphics[scale=0.15]{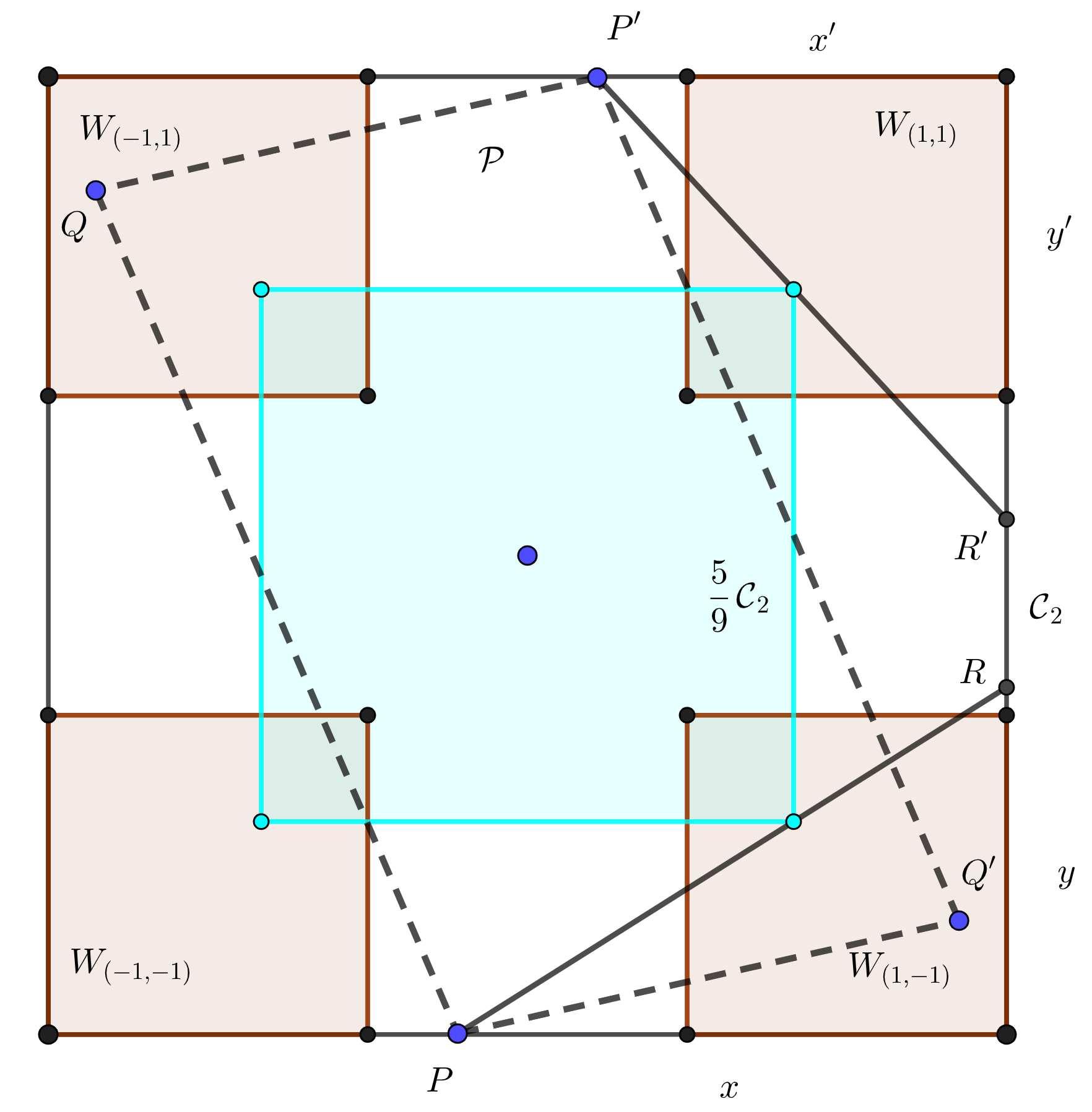}
    \caption{Proof of Lemma \ref{lemkwadrat}. Point $R$ is the intersection of a line passing through $P$ and $\frac{5}{9}(1, -1)$ with a side of $\mathcal{C}_2$. Similarly  $R'$ is the intersection of a line passing through $P'$ and $\frac{5}{9}(1, 1)$ with the same side of $\mathcal{C}_2$. If $y$ and $y'$ are distances of $R$ to $(1, -1)$ and $R'$ to $(1, 1)$ respectively, then by a direct calculation we obtain an inequality $y+y'<2$. This is a contradiction -- for the parallelogram $\mathcal{P}$ to contain $\frac{5}{9}\mathcal{C}_2$ it would be necessary that $Q'$ lies below the line $PR$ and above the line $P'R'$ at the same time. This is impossible, as these two lines intersect outside the square $\mathcal{C}_2$.} 
    \label{ryskwadrat3}
\end{figure}

Let $x$ be the length of the segment connecting $P$ with the vertex $(1, -1)$ of $\mathcal{C}_2$ and $x'$ be the length of the segment connecting $P'$ with the vertex $(1, 1)$ (see Figure \ref{ryskwadrat3}). By the assumption we have $x, x' \in \left ( \frac{2}{3}, \frac{4}{3} \right )$ and also $x+x'=2$, since $P$ and $P'$ are $0$-symmetric. Let $R$ be the point of intersection of the line passing through $P$ and $\frac{5}{9}(1, -1)$ with the side $[(1, 1), (1, -1)]$ of $\mathcal{C}_2$, and similarly let $R'$ be the point of intersection of the line passing trough $P'$ and $\frac{5}{9}(1, 1)$ with the same side of $\mathcal{C}_2$. By $y$ and $y'$ we denote the distances of $R$ to $(1, -1)$ and $R'$ to $(1, 1)$ respectively. By a simple calculation we get
$$y = \frac{4x}{9x-4} \quad \text{ and } \quad y' = \frac{4x'}{9x'-4}.$$
Hence
$$y + y' = \frac{8}{9} + \frac{16}{9} \cdot \frac{10}{(9x-4)(14x-9)}.$$
However, for $x \in \left ( \frac{2}{3}, \frac{4}{3} \right )$ we have $(9x-4)(14x-9)>16.$ Thus
$$y + y' < \frac{8}{9}  + \frac{10}{9} = 2.$$
This proves that it is impossible to complete the points $P$ and $P'$ to a $0$-symmetric parallelogram containing $\frac{5}{9}\mathcal{C}_2$. This contradicts our assumption and the conclusion follows.

\begin{figure}
    \centering
    \includegraphics[scale=0.15]{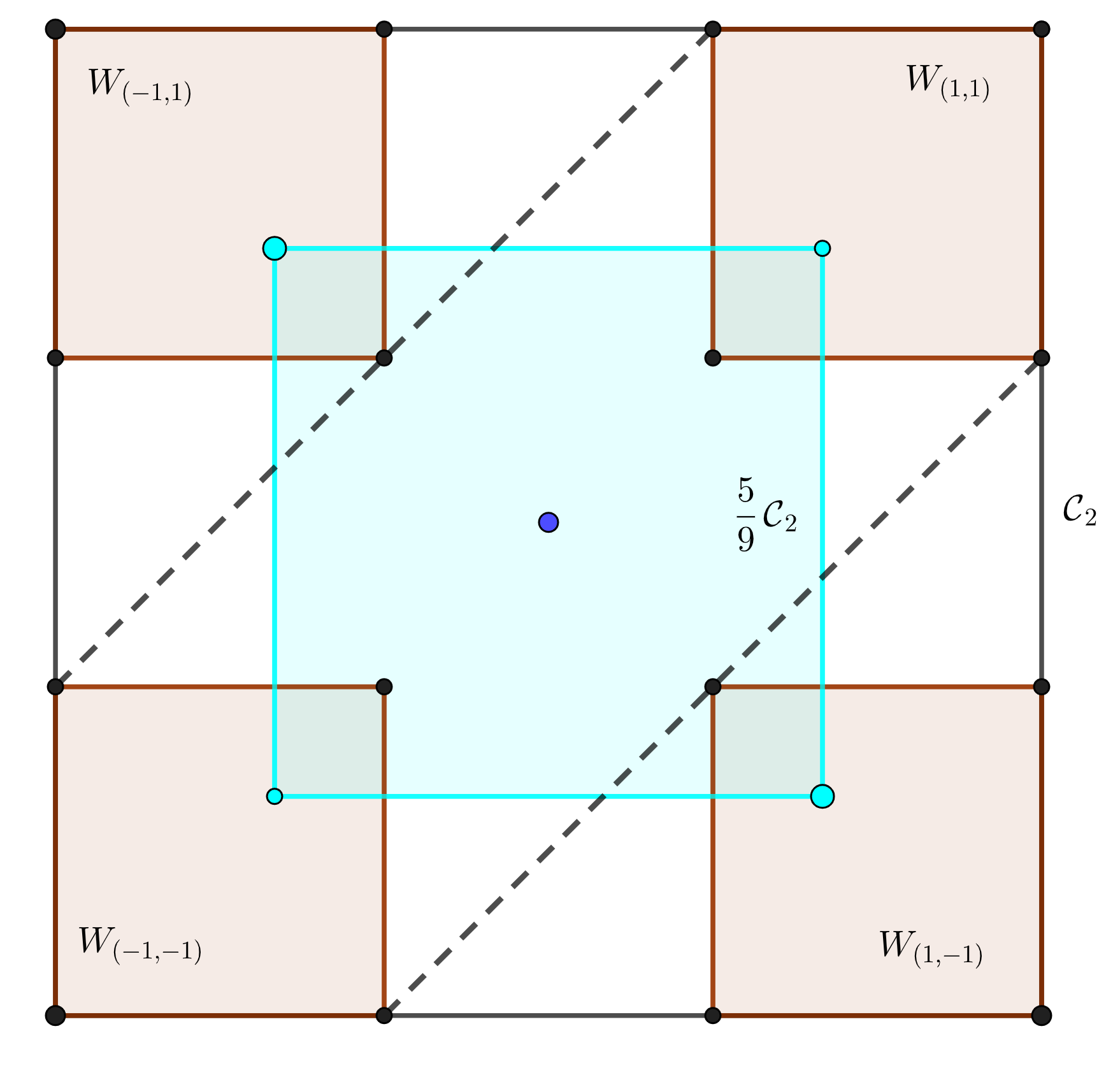}
    \caption{Proof of Lemma \ref{lemkwadrat}. If the squares $W_{(1, -1)}$ and $W_{(-1, 1)}$ do not contain any vertices of the parallelogram $\mathcal{P}$, then $\mathcal{P}$ has to be contained in the region bounded by the two dashed lines. In this case, the vertices $\frac{5}{9}(-1, 1)$ and $\frac{5}{9}(1, -1)$ of $\frac{5}{9}\mathcal{C}_2$ are not in $\mathcal{P}$.}
    \label{ryskwadrat4}
\end{figure}

We are left with proving that each of the squares $W_{(\varepsilon_1, \varepsilon_2)}$ is non-empty. Let us assume the opposite. Because we have just proved that each vertex of $\mathcal{P}$ belongs to some square $W_{(\varepsilon_1, \varepsilon_2)}$, two $0$-symmetric squares $W_{(\varepsilon_1, \varepsilon_2)}$ have to contain two vertices of $\mathcal{P}$ each. Hence, the parallelogram $\mathcal{P}$ is contained in the region bounded by two dashed lines, as presented in Figure \ref{ryskwadrat4} (or in the analogous region along the other diagonal of $\mathcal{C}_2$).

In this case however, the vertices $\frac{5}{9}(-1, 1)$ and $\frac{5}{9}(1, -1)$ of the square $\frac{5}{9}\mathcal{C}_2$ are outside of $\mathcal{P}$ and the assumed inclusion $\frac{5}{9}\mathcal{C}_2 \subseteq \mathcal{P}$ does not hold. We have again obtained a contradiction and the proof is finished. \qed

We are ready to prove the main result of this section.

\begin{twr}
\label{twr3d}
We have the equality $d_{BM}(\mathcal{C}_3, \mathcal{C}_3^*) = \frac{9}{5}$. Moreover, if a linear operator $T: \mathbb{R}^3 \to \mathbb{R}^3$ satisfies $\frac{5}{9} \mathcal{C}_3 \subseteq T(\mathcal{C}_3^*) \subseteq \mathcal{C}_3$, then the matrix of $T$ is of the form
$$\begin{bmatrix}
\frac{1}{3} & -1 & -1\\
-1 & \frac{1}{3} & -1\\
-1 & -1 & \frac{1}{3} \\
\end{bmatrix}$$
or arises from the matrix above by operations of: permuting of rows/columns, multiplying a row/column by $-1$ (there are in total $192$ of such matrices). 
\end{twr}

\emph{Proof.} We will say that an octahedron $K \subseteq \mathbb{R}^3$ is \emph{nice} if there exists a vertex $a$ of the unit cube $\mathcal{C}_3$, such that the vertices of $K$ are of the form $ \pm \left ( \frac{1}{3}a + \frac{2}{3}b_i \right )$, where $i=1,2,3$ and $b_i$ are the vertices of the cube adjacent to $a$.  We note that there exist a total of $4$ nice octahedrons since the vertex $a$ can be chosen in $8$ ways and two symmetric choices of $a$ give rise to the same nice octahedron. 

Our goal is to prove that for every linear operator $T: \mathbb{R}^3 \to \mathbb{R}^3$ the following implication holds
\begin{equation}
\frac{5}{9}\mathcal{C}_3 \subseteq T(\mathcal{C}_3^*) \subseteq \mathcal{C}_3 \implies T(\mathcal{C}_3^*) \text{ is nice}.
\end{equation}

This is sufficient for establishing the equality $d_{BM}(\mathcal{C}_3, \mathcal{C}_3^*) = \frac{9}{5}$, as for a nice octahedron $T(\mathcal{C}_3^*)$ it is straightforward to verify that $r\mathcal{C}_3 \nsubseteq T(\mathcal{C}_3^*)$ for $r > \frac{5}{9}$. Moreover, it can be easily checked that the octahedron $T(\mathcal{C}_3^*)$ is nice if and only if $T$ has the matrix representation indicated in the statement. From here on we will assume that the linear operator $T:\mathbb{R}^3 \to \mathbb{R}^3$ satisfies the inclusions $\frac{5}{9}\mathcal{C}_3 \subseteq T(\mathcal{C}_3^*) \subseteq \mathcal{C}_3$.

Let $V \subseteq \mathbb{R}^3$ be the set defined as
$$V = \left \{ (x,y,z) \in \mathbb{R}^3 : \frac{1}{3} \leq |x|, |y|, |z| \leq 1  \right \}.$$
The set $V$ is a union of $8$ disjoint closed cubes in $\mathbb{R}^3$,  each containing a unique vertex of $\mathcal{C}_3$. For $\varepsilon_1, \varepsilon_2, \varepsilon_3 \in \{-1, 1\}$ let $V_{(\varepsilon_1, \varepsilon_2, \varepsilon_3)}$ be the cube containing the vertex $(\varepsilon_1, \varepsilon_2, \varepsilon_3).$ In other words
$$V_{(\varepsilon_1, \varepsilon_2, \varepsilon_3)}= \left \{ (x,y,z) \in \mathbb{R}^3 : \frac{1}{3} \leq \varepsilon_1x, \varepsilon_2y, \varepsilon_3z \leq 1 \right \}.$$
We will prove the following claim, closely resembling Lemma \ref{lemkwadrat} in the three-dimensional setting.

\textbf{Claim}. Each vertex of $T(\mathcal{C}_3^*)$ belongs to $V$. Moreover, every cube $V_{(\epsilon_1, \epsilon_2, \epsilon_3)}$ (where $\epsilon_{i} \in \{-1, 1\}$) contains at most one vertex of $T(\mathcal{C}_3^*)$.

Let $e_1, e_2, e_3$ be the standard basis in $\mathbb{R}^3$. Let $H \subseteq \mathbb{R}^3$ be a two-dimensional subspace of $\mathbb{R}^3$ (a plane passing through $0$) such that $e_i \in H$ for some $1 \leq i \leq 3$.

To prove our Claim, we will rely on the following straightforward observation: the orthogonal projection of $\mathcal{C}_3$ onto $H$ is a rectangle, the orthogonal projection of $\frac{5}{9}\mathcal{C}_3$ is the same rectangle scaled by $\frac{5}{9}$ and the $i$-th coordinate is preserved by this projection.

Let $w_1, w_2, w_3$ be different, non-symmetric vertices of $T(\mathcal{C}_3^{*})$ such that $0 \leq z_1 \leq z_2 \leq z_3$ where $w_j = (x_j, y_j, z_j)$ for $j=1,2,3$. We take a vector $f$ perpendicular to the plane through $w_1, w_2, w_3$ and $H \subseteq \mathbb{R}^3$ as a two-dimensional subspace containing $e_3$ and $f$. If by $P: \mathbb{R}^3 \to H$ we denote the orthogonal projection to $H$, then the plane through $w_1, w_2, w_3$ is projected onto a line. Since $0\leq z_1 \leq z_2 \leq z_3$ and $z$ coordinate is preserved under projection, $P(w_2)$ is between $P(w_1)$ and $P(w_3)$ on this line, so that $P(w_2) \in [P(w_1), P(w_3)]$. Because $T(\mathcal{C}_3^*) = \conv \{ \pm w_1, \pm w_2, \pm w_3 \}$, we conclude that $P(T(\mathcal{C}_3^*)) = \conv \{ \pm P(w_1), \pm P(w_3) \}$. In other words, the projection of the octahedron $T(\mathcal{C}_3^*)$ is a parallelogram with the vertices $\pm P(w_1), \pm P(w_3)$. We also have $\frac{5}{9} P(\mathcal{C}_3) \subseteq P(T(\mathcal{C}_3^*)) \subseteq P(\mathcal{C}_3)$ so by Lemma \ref{lemkwadrat}, we conclude that $z_1, z_3 \geq \frac{1}{3}$, and thus $z_2 \geq \frac{1}{3}$ as well. Note that $P(\mathcal{C}_3)$ is a rectangle and not necessarily a square. However, Lemma \ref{lemkwadrat} can still be applied, by first transforming $P(\mathcal{C}_3)$ into a square while preserving the $z$ coordinate, which is done by scaling along the other coordinate axis in the projection plane.  An analogous reasoning, applied for the other two coordinates $x$ and $y$, yields the first part of our Claim.

To finish the proof of the Claim, we are left with showing that no two vertices of $T(\mathcal{C}_3^*)$ are in the same cube $V_{(\epsilon_1, \epsilon_2, \epsilon_3)}$. Assuming the opposite, two possibilities emerge. Either three non-symmetric vertices of $T(\mathcal{C}_3^*)$ are all in the same cube, or the three vertices are contained in two cubes, which can be separated from their symmetric copy by a plane parallel to some face of the cube $\mathcal{C}_3$ (spanned by $e_i, e_j$ for some $i \neq j$). In the latter case, we take $H'=\lin \{e_i, e_j \}$, and in the former we take any $i \neq j$ and $H'$ defined in the same way. If now $P': \mathbb{R}^3 \to \mathbb{R}^3$ is the orthogonal projection onto $H'$, then again we have $\frac{5}{9} P'(\mathcal{C}_3) \subseteq P'(T(\mathcal{C}_3^*)) \subseteq P'(\mathcal{C}_3)$, but all vertices of $T(\mathcal{C}_3^*)$ will be projected onto two opposite squares (that are defined like in the statement of Lemma \ref{lemkwadrat}). In this case we can not refer to Lemma \ref{lemkwadrat} directly, as we do not know if the projection of $T(\mathcal{C}_3^*)$ is a parallelogram. However, we obtain a contradiction in exactly the same way as in the last step of the proof of this lemma (see Figure \ref{ryskwadrat4}).

With our Claim proved, we are ready to finish the proof. As there are $8$ cubes in $V$ and $6$ vertices of $T(\mathcal{C}_3^*)$, there exists a cube in $V$ not containing any vertex of $T(\mathcal{C}_3^*)$. Without loss of generality, we can assume that the cubes $V_{(-1, -1, -1)}$, $V_{(1, 1, 1)}$ do not contain any vertex of $T(\mathcal{C}_3^*)$.

The other six cubes in $V$ contain some vertex of $T(\mathcal{C}_3^*)$, so now we let $v_1, v_2, v_3$ be the vertices of $T(\mathcal{C}_3^*)$ in $V_{(1, -1, -1)}, V_{(-1, 1, -1)}, V_{(-1, -1, 1)}$ respectively. Let us also write $v_i = (x_i, y_i, z_i)$ for $1 \leq i \leq 3$ (here we do not assume that $z_i$ are non-negative or ordered, as previously in the proof of the Claim). From the condition $\frac{5}{9} \mathcal{C}_3 \subseteq T(\mathcal{C}_3^*)$ it follows now that $\left \| \frac{v_1+v_2+v_3}{3} \right \|_{\infty} \geq \frac{5}{9}$. Still without losing the generality we can suppose that
$$\left | \frac{z_1+z_2+z_3}{3} \right | = \left \| \frac{v_1+v_2+v_3}{3} \right \|_{\infty} \geq \frac{5}{9}.$$
Since $v_3 \in \mathcal{C}_3$ we have $z_3 \leq 1$, and since $v_1 \in V_{(1, -1, -1)}$ and $v_2 \in V_{(-1, 1, -1)}$ we have
$$\frac{z_1+z_2+z_3}{3} \leq \frac{1}{3} \left ( - \frac{1}{3} - \frac{1}{3} + 1 \right )  = \frac{1}{9} < \frac{5}{9}.$$
Thus
$$\frac{5}{9} \leq - \frac{z_1+z_2+z_3}{3} \leq -\frac{1}{3} \left ( -1 -1 + \frac{1}{3} \right )=\frac{5}{9}.$$
Hence we have an equality, and it follows that $z_1=z_2=-1$ and $z_3 = \frac{1}{3}$.

For sufficiently small $\varepsilon>0$ we have 
$$\left ( \frac{1}{3} - \varepsilon \right ) v_1 + \left ( \frac{1}{3} - \varepsilon \right ) v_2 + \left ( \frac{1}{3} + 2\varepsilon \right ) v_3 \in T(\mathcal{C}_3^*)$$ 
and thus
$$\left \|\left ( \frac{1}{3} - \varepsilon \right ) v_1 + \left ( \frac{1}{3} - \varepsilon \right ) v_2 + \left ( \frac{1}{3} + 2\varepsilon \right ) v_3 \right \|_{\infty} \geq \frac{5}{9}.$$
However, if $\varepsilon>0$ is small enough, then
$$\left | \left ( \frac{1}{3} - \varepsilon \right ) z_1 + \left ( \frac{1}{3} - \varepsilon \right ) z_2 + \left ( \frac{1}{3} + 2\varepsilon \right ) z_3 \right |= -  \left (\left ( \frac{1}{3} - \varepsilon \right ) z_1 + \left ( \frac{1}{3} - \varepsilon \right ) z_2 + \left ( \frac{1}{3} + 2\varepsilon \right ) z_3 \right )$$
$$=\left ( \frac{1}{3} - \varepsilon \right ) + \left ( \frac{1}{3} - \varepsilon \right ) - \left (\frac{1}{9} + \frac{2}{3} \varepsilon \right ) = \frac{5}{9} - \frac{4}{3}\varepsilon < \frac{5}{9}.$$
Thus, the maximum norm $\| \left ( \frac{1}{3} - \varepsilon \right ) v_1 + \left ( \frac{1}{3} - \varepsilon \right ) v_2 + \left ( \frac{1}{3} + 2\varepsilon \right ) v_3 \|_\infty$ has to be attained on one of the other two coordinates. By taking $\varepsilon = \frac{1}{N}$ and letting $N \to \infty$, we see that one of the coordinates realizes the maximum infinitely many times. We can suppose that it is the $y$ coordinate. Thus, for infinitely many $N \geq 1$ we have
$$\left | \left ( \frac{1}{3} - \frac{1}{N} \right ) y_1 + \left ( \frac{1}{3} - \frac{1}{N} \right ) y_2 + \left ( \frac{1}{3} + \frac{2}{N} \right ) y_3 \right | \geq \frac{5}{9}.$$
By taking $N \to \infty$ and passing to the limit, we get
$$\left | \frac{y_1+y_2+y_3}{3} \right |\geq \frac{5}{9}.$$
Reasoning exactly like before we prove that $\frac{y_1+y_2+y_3}{3} \leq \frac{1}{9}$, implying $\frac{y_1+y_2+y_3}{3} \leq -\frac{5}{9}$ and then we obtain $y_1=y_3=-1$, $y_2=\frac{1}{3}$.

To finish the proof, we observe that now for sufficiently small $\varepsilon>0$, we have again that $\left ( \frac{1}{3} - 2\varepsilon \right ) v_1 + \left ( \frac{1}{3} + \varepsilon \right ) v_2 + \left ( \frac{1}{3} +\varepsilon \right ) v_3 \in T(\mathcal{C}_3^*)$ and
$$\left |\left ( \frac{1}{3} - 2\varepsilon \right ) z_1 + \left ( \frac{1}{3} + \varepsilon \right ) z_2 + \left ( \frac{1}{3} +\varepsilon \right ) z_3 \right | = - \left ( \left ( \frac{1}{3} - 2\varepsilon \right ) z_1 + \left ( \frac{1}{3} + \varepsilon \right ) z_2 + \left ( \frac{1}{3} +\varepsilon \right ) z_3 \right )$$
$$= \left ( \frac{1}{3} - 2\varepsilon \right ) + \left ( \frac{1}{3} + \varepsilon \right ) - \left ( \frac{1}{9} + \frac{1}{3}\varepsilon \right ) = \frac{5}{9}  - \frac{4}{9} \varepsilon<\frac{5}{9}.$$
In exactly the same way we estimate the absolute value of the $y$ coordinate. This shows that the maximum norm has to be achieved on the $x$ coordinate and thus
$$\left | \frac{x_1+x_2+x_3}{3} \right |\geq \frac{5}{9}.$$
The same argument as before now gives us that $x_1 = \frac{1}{3}$ and $x_2=x_3=-1$. This shows that the octahedron $T(\mathcal{C}_3^*)$ is nice and the conclusion follows. \qed
\section{Banach-Mazur distance between the cube and the cross-polytope in the four-dimensional case}
\label{sect4d}

In this section, we prove that $d_{BM}(\mathcal{C}_4, \mathcal{C}_4^*)=2$, again confirming a conjecture suggested by numerical data. Interestingly, the proof in dimension four does not seem to share much similarity with the three-dimensional case. 

In dimension four we do not characterize all operators $T$ such that $\frac{1}{2}\mathcal{C}_4 \subseteq T(\mathcal{C}_4^*) \subseteq \mathcal{C}_4$, but we provide some examples:

$$\begin{bmatrix}
    1 & 1 & 1 & 0 \\
    1 & 1 & -1 & 0 \\
    1 & -1 & 0 & 1 \\
    1 & -1 & 0 & -1 
    \end{bmatrix}, \quad
    \begin{bmatrix}
    1 & 1 & 1 & 0 \\
    1 & -1 & 0 & 1 \\
    1 & 0 & -1 & -1 \\
    0 & 1 & -1 & 1
    \end{bmatrix}, \quad
    \begin{bmatrix}
    1 & 1 & 1 & -1 \\
    -1 & 1 & 1 & 1 \\
    1 & -1 & 1 & 1 \\
    1 & 1 & -1 & 1
    \end{bmatrix}.
    $$
These three linear operators yield the upper bound of $2$ and are essentially different from each other. This contrasts with the three-dimensional setting.

It should be noted that in the proof we will use a rather unusual meaning of $\sgn(x)$. We define $\sgn(x)=1$ for $x \geq 0$ and $\sgn(x)=-1$ for $x<0$. Thus $\sgn(x) \in \{1, -1\}$ for every $x \in \mathbb{R}$ and $\sgn(x) x = |x|$.

\begin{twr}
\label{twr4d}
We have the equality $d_{BM}(\mathcal{C}_4, \mathcal{C}_4^*) = 2$.
\end{twr}

\emph{Proof.} We have already mentioned some examples of operators $T$ providing the upper bound $d_{BM}(\mathcal{C}_4, \mathcal{C}_4^*) \leq 2$, so our goal is to establish the opposite estimate. With the aim of obtaining a contradiction, we assume that $d_{BM}(\mathcal{C}_4, \mathcal{C}_4^*) < 2$. This means that there exists a linear operator $T: \mathbb{R}^4 \to \mathbb{R}^4$ such that $r\mathcal{C}_4 \subseteq T(\mathcal{C}_4^*) \subseteq \mathcal{C}_4$ where $r>\frac{1}{2}$. From the fact that $\|T(e_i)\|_{\infty} \leq 1$ for $1 \leq i \leq 4$, it follows that the absolute value of each entry of the matrix associated with $T$ is not greater than $1$. We will denote the rows of this matrix as $x, y, z, w \in \mathbb{R}^4$ (where $x=(x_1, x_2, x_3, x_4)$ and similarly for the other rows). We can perform several operations on the operator $T$ satisfying $r\mathcal{C}_4 \subseteq T(\mathcal{C}_4^*) \subseteq \mathcal{C}_4$ that preserve the inclusions. These include: swapping rows, swapping columns, changing the sign of all elements in a column, changing the sign of all elements in a row. Moreover, the inclusions are also preserved when a column or a row is multiplied by a real number $\lambda>1$, assuming that after the multiplication each entry of the matrix of $T$ still has absolute value not greater than $1$. Thus, we can assume the following:
\begin{equation}
\label{propkolumna}
\text{Each column of $T$ contains an element with the absolute value equal to } 1.
\end{equation}
We shall call a vector $s = (s_1, s_2, s_3, s_4) \in \mathbb{R}^4$ a \textit{string}, if it lies on the boundary of $\mathcal{C}_4^*$, or in other words if it satisfies the condition $|s_1| + |s_2| + |s_3| + |s_4| = 1$. We will refer shortly to a pair of symmetric strings $(s, -s)$ as a \textit{string pair}.

If $s$ is a string, then the point $T(s)$ is on the boundary of $T(\mathcal{C}_4^*)$ which implies that $\|T(s)\|_\infty \geq r$ or equivalently, there exists a row $a \in \{x, y, z, w \}$ such that $| \langle s, a \rangle | \geq r$. In this case, we will say that the string $s$ is \emph{associated} to the row $a$, and we will write shortly $a \sim s$. Hence, each string is associated with at least one row. Clearly, if $a \sim s$, then also $a \sim (-s)$, so we can speak of string pairs associated with a given row.

Now, let $\mathcal{S} = \left \{ (\pm \frac{1}{4}, \pm \frac{1}{4}, \pm \frac{1}{4}, \pm \frac{1}{4}) \right \}$ be a set of $16$ strings. We will establish the following properties:
\begin{equation}
\text{Each string in the set $\mathcal{S}$ is associated with exactly one row.}
\end{equation}
\begin{equation}
\text{Each row is associated with exactly $4$ strings from $\mathcal{S}$}.
\end{equation}
 To prove the properties above, let us assume that for a row $a$ we have $a = (a_1, a_2, a_3, a_4)$ and $a_1 \geq a_2 \geq a_3 \geq a_4 \geq 0$. If $a$ is associated with at least $5$ strings from $\mathcal{S}$, then it is associated with at least $3$ string pairs. Clearly, among these pairs there are the following string pairs: $ \pm (\frac{1}{4}, \frac{1}{4}, \frac{1}{4}, \frac{1}{4}) $, $\pm (\frac{1}{4}, \frac{1}{4}, \frac{1}{4}, - \frac{1}{4})$, and $\pm (\frac{1}{4}, \frac{1}{4},- \frac{1}{4}, \frac{1}{4})$ as these maximize the value $| \langle s, a \rangle |$. However, looking at the last string pair, we get
 $$0 \leq \left \langle \left ( \frac{1}{4}, \frac{1}{4},- \frac{1}{4}, \frac{1}{4} \right ), a \right  \rangle = \frac{1}{4} (a_1+a_2-a_3+a_4) \leq \frac{1}{4} (a_1+a_2) \leq \frac{1}{4} \cdot 2 = \frac{1}{2} < r,$$
 which gives us a contradiction. Thus we have proved that each row has at most $2$ string pairs from $\mathcal{S}$ associated, or in other words, at most $4$ strings. Because there are $4$ rows and $16$ strings in $\mathcal{S}$, it follows from simple counting that each row is associated with exactly $4$ strings and each string with exactly one row (as each string has to be associated with some row). Here it should be noted that during the latter part of the reasoning, we will often refer to the fact, that for a given row we know exactly the two string pairs from $\mathcal{S}$ associated to it: again, if $a = (a_1, a_2, a_3, a_4)$ is a row and $a_1 \geq a_2 \geq a_3 \geq a_4 \geq 0$, then these are the string pairs: $ \pm (\frac{1}{4}, \frac{1}{4}, \frac{1}{4}, \frac{1}{4}) $, $\pm (\frac{1}{4}, \frac{1}{4}, \frac{1}{4}, - \frac{1}{4})$. Since the string $(\frac{1}{4}, \frac{1}{4}, - \frac{1}{4}, \frac{1}{4})$ can not be associated to $a$ in this case, we must have that $a_3>a_4$. This implies the next important observation:
 \begin{equation}
     \text{In every row, there exists a \emph{unique} element with the minimal absolute value.}
 \end{equation}
 
 Now we will take a closer look at the entries of the matrix representing $T$. To visualize the possible situations more clearly, we will make the following identifications:
    \begin{itemize}
        \item For each row, we will denote the unique element with the minimal absolute value as $\circ$. Such an element will be called minimal.
        \item For each row, the non-minimal elements are non-zero, so we will denote the positive elements as $+$ and the negative elements as $-$.
    \end{itemize}
In the beginning of the proof, we have mentioned several operations that can be done on the matrix representing $T$, preserving the inclusions $r\mathcal{C}_4 \subseteq T(\mathcal{C}_4^*) \subseteq \mathcal{C}_4$. We shall classify all possible matrices of $T$ under the just defined identification and the mentioned symmetries. In order to do this, we shall use the following two rules for any two different rows $a,b \in \{x, y, z, w\}$.
\begin{enumerate}[(i)]
    \item \label{rule1} The rows $a,b$ do not have a matching layout in terms of signs.
    \item \label{rule2} If rows $a,b$ have minimal elements in two different columns, then they do not match in the remaining two columns.
\end{enumerate}

To see that~\eqref{rule1} is true, note that we can determine the $4$ associated strings from the sign layout of a row. Thus if two rows have the same layout, they would be associated to the same strings which is a contradiction as each string is associated to exactly one row. To show that~\eqref{rule2} holds, suppose that there are two rows $a,b$ with signs laid out as $\begin{bmatrix}
    + & + &  & \circ \\
    + & + & \circ  & 
    \end{bmatrix}.$ In this case, we would have two rows associated to the same string $\frac{1}{4}(1, 1, \sgn(a_3), \sgn(b_4))$, which is a contradiction. We can apply similar reasoning for every possible layout of signs in the first two columns, thus proving~\eqref{rule2}. 
    
We consider the column of $T$ with the largest number of minimal elements. Without loss of generality, we can assume that it is the last (fourth) column. There are four cases possible -- the last column containing exactly $i$ minimal elements, where $1 \leq i \leq 4$. We will denote the application of rules as $\xrightarrow{(\ast)[a,b]}$ where $\ast$ is the rule number (\ref{rule1}, or \ref{rule2}) and $a,b$ are rows to which the rule is applied to.

\begin{enumerate}
    \item \textbf{$4$ minimal elements}. This case is easy to discard, as the last column has all elements with the absolute value smaller than $1$, which contradicts the assumption (\ref{propkolumna}).
    \item \textbf{$3$ minimal elements}. In this case, by using the aforementioned symmetries, the matrix of $T$ can be represented as the leftmost matrix below. For example, we can first move the minimal elements to the desired place by permuting rows/columns, and then adequately adjust the signs by multiplying rows/columns with $-1$. Then we apply the two rules stated previously to determine other entries.
    $$\begin{bmatrix}
    + & + & + & \circ \\
    + &  &  & \circ \\
    + &  &  & \circ \\
    + &  & \circ  & 
    \end{bmatrix} \xrightarrow{\eqref{rule2}[x,w]}
    \begin{bmatrix}
    + & + & + & \circ \\
    + &  &  & \circ \\
    + &  &  & \circ \\
    + & -  & \circ  & 
    \end{bmatrix} \xrightarrow{\substack{\eqref{rule2}[y,w] \\ \eqref{rule2}[z,w]}}
\begin{bmatrix}
    + & + & + & \circ \\
    + & + &  & \circ \\
    + & + &  & \circ \\
    + & -  & \circ  & 
    \end{bmatrix} \xrightarrow{\substack{\eqref{rule1}[x,y] \\ \eqref{rule1}[x,z]}}
    \begin{bmatrix}
    + & + & + & \circ \\
    + & + & - & \circ \\
    + & + & - & \circ \\
    + & -  & \circ  & 
    \end{bmatrix}
$$
We arrive at a contradiction as the rule~\eqref{rule1} is violated for rows $y, z$. Thus, we can discard this case.

    \item \textbf{$2$ minimal elements}. This case can be subdivided into further two essentially different possibilities -- if all the minimal elements are contained in two or three different columns. We start with the latter case. In this situation, we can assume that the matrix of $T$ has the form as the leftmost matrix.
    $$\begin{bmatrix}
    + & + & + & \circ \\
    + &  &  & \circ \\
    + &  & \circ  &  \\
    + & \circ  &   & 
    \end{bmatrix} \xrightarrow{\substack{\eqref{rule2}[x,z] \\ \eqref{rule2}[x,w]}}
    \begin{bmatrix}
    + & + & + & \circ \\
    + &  &  & \circ \\
    + & - & \circ  &  \\
    + & \circ  & -  & 
    \end{bmatrix} \xrightarrow{\substack{\eqref{rule2}[y,z] \\ \eqref{rule2}[y,w]}}
    \begin{bmatrix}
    + & + & + & \circ \\
    + & + & + & \circ \\
    + & - & \circ  &  \\
    + & \circ  & -  & 
    \end{bmatrix}
    $$
    Since rule~\eqref{rule1} is violated for rows $x,y$ we discard this possibility.
    
    If all the minimal elements are contained in two different columns, then we can assume the following form of the matrix of $T$:
    $$\begin{bmatrix}
    + & + & + & \circ \\
    + &  &  & \circ \\
    + &  & \circ  & +  \\
    + &  & \circ  & 
    \end{bmatrix} \xrightarrow{\substack{\eqref{rule2}[x, z] \\ \eqref{rule2}[x, w]}}
    \begin{bmatrix}
    + & + & + & \circ \\
    + &  &  & \circ \\
    + & - & \circ  &  + \\
    + & - & \circ   & 
    \end{bmatrix} \xrightarrow{\substack{\eqref{rule1}[z,w] \\ \eqref{rule2}[y,z]}}
    \begin{bmatrix}
    + & + & + & \circ \\
    + & + &  & \circ \\
    + & - & \circ  &  + \\
    + & - & \circ   & -
    \end{bmatrix} \xrightarrow{\eqref{rule1}[x,y]} 
    \begin{bmatrix}
    + & + & + & \circ \\
    + & + & - & \circ \\
    + & - & \circ  &  + \\
    + & - & \circ   & -
    \end{bmatrix}.
$$

The final matrix satisfies both rules, so this situation will require a further examination which will be carried out in the next part of the proof.
\item \textbf{$1$ minimal element}. In this case, the matrix of $T$ can be brought into the following form:

$$
\begin{bmatrix}
    + & + & + & \circ \\
    + &  & \circ  & + \\
    + & \circ  &  &   \\
    \circ & +  &    & 
    \end{bmatrix} \xrightarrow{\substack{\eqref{rule2}[x,y] \\ \eqref{rule2}[x,z] \\ \eqref{rule2}[x,w]}}
    \begin{bmatrix}
    + & + & + & \circ \\
    + & - & \circ  & + \\
    + & \circ  & - &   \\
    \circ & +  &  -  &  
    \end{bmatrix}\xrightarrow{\substack{\eqref{rule2}[y,z]}}
    \begin{bmatrix}
    + & + & + & \circ \\
    + & - & \circ  & + \\
    + & \circ  & - &  - \\
    \circ & +  &  -  &  
    \end{bmatrix}\xrightarrow{\substack{\eqref{rule2}[z,w]}}
    \begin{bmatrix}
    + & + & + & \circ \\
    + & - & \circ  & + \\
    + & \circ  & - &  - \\
    \circ & +  &  -  &  +
    \end{bmatrix}.
$$
This matrix also satisfies both rules and also needs to be examined further.
\end{enumerate}

We should recall here that for $r = \frac{1}{2}$ the two matrices
$$\begin{bmatrix}
    1 & 1 & 1 & 0 \\
    1 & 1 & -1 & 0 \\
    1 & -1 & 0 & 1 \\
    1 & -1 & 0 & -1 
    \end{bmatrix} \quad \text{ and } \quad
    \begin{bmatrix}
    1 & 1 & 1 & 0 \\
    1 & -1 & 0 & 1 \\
    1 & 0 & -1 & -1 \\
    0 & 1 & -1 & 1
    \end{bmatrix}
    $$
do satisfy the desired inclusions and they have exactly the respective form. This somewhat explains why these two situations require more work to be done. We remind that we have assumed that $r>\frac{1}{2}$ and we are aiming at the contradiction.

We start with the second possibility, that is when the matrix of $T$ is of the form
 $$   \begin{bmatrix}
    + & + & + & \circ \\
    + & -  & \circ  & + \\
    + & \circ  & -  &  -  \\
    \circ & +  & -   & + 
    \end{bmatrix}.
    $$
For each row we shall construct a specific string associated to it. We start with the row $x$ and we define a function $s_x: [0, 1] \to \bd \mathcal{C}_4^*$ as:
$$s_x(t) = \left ( \frac{1-t}{3}, \frac{1-t}{3}, \frac{1-t}{3}, -\sgn(x_4)t \right )$$
for $t \in [0, 1]$. By the definition, the vector $s_x(t)$ is a string for every $t \in [0, 1]$. We note that
$$\langle x, s_x(t) \rangle = \frac{1-t}{3}x_1 + \frac{1-t}{3}x_2 + \frac{1-t}{3}x_3 - t|x_4|,$$
where $\langle \cdot, \cdot \rangle$ denotes the standard scalar product in $\mathbb{R}^4$. Thus for $t = \frac{1}{4}$ we get
$$\left \langle x, s_x\left ( \frac{1}{4} \right ) \right \rangle = \frac{1}{4} x_1 + \frac{1}{4} x_2 + \frac{1}{4} x_3 - \frac{1}{4} |x_4| \geq r$$
since $x \sim (\frac{1}{4}, \frac{1}{4}, \frac{1}{4}, \pm \frac{1}{4}) \in \mathcal{S}$. On the other hand, for $t = \frac{1}{2}$ the opposite inequality is true, as
$$0 \leq\left \langle x, s_x\left ( \frac{1}{2} \right ) \right \rangle = \frac{1}{6}x_1 + \frac{1}{6}x_2 + \frac{1}{6}x_3 - \frac{1}{2}|x_4| \leq 3 \cdot \frac{1}{6} = \frac{1}{2} < r.$$
Because the function $\langle x, s_x(t) \rangle$ is linear, there exists a unique $t_x \in [ \frac{1}{4}, \frac{1}{2})$ such that $\langle x, s_x(t_x) \rangle = r$. Moreover we have $\langle x, s_x(t) \rangle \geq r \Leftrightarrow t \leq t_x$. We will call the vector $s_x(t_x)$ the \emph{specific string} of the row $x$. The functions $s_y, s_z, s_w: [0, 1] \to \bd \mathcal{C}_4^*$ are defined similarly. More precisely:
$$s_y(t) = \left ( \frac{1-t}{3}, -\frac{1-t}{3}, -\sgn(y_3)t, \frac{1-t}{3} \right ),$$
$$s_z(t) = \left ( \frac{1-t}{3}, -\sgn(z_2)t, -\frac{1-t}{3}, -\frac{1-t}{3} \right ),$$
$$s_w(t) = \left ( -\sgn(w_1)t, \frac{1-t}{3}, -\frac{1-t}{3}, \frac{1-t}{3} \right ). $$
The numbers $t_y, t_z, t_w$ are unique numbers in $[ \frac{1}{4}, \frac{1}{2})$ such that
$$s_y(t_y)=s_z(t_z)=s_w(t_w)=r$$
and $s_y(t_y), s_z(t_z), s_w(t_w)$ are called the specific strings of rows $y, z, w$ respectively.

By definition, for every row $a \in \{x, y, z, w\}$ the specific string $s_a(t_a)$ is associated to $a$. The crucial property of specific strings, that we are going to establish, is the following: for every row $a$, the specific string $s_a(t_a)$ of $a$ is also associated to some other row other than $a$. Indeed, let $a$ be a fixed row. Then for sufficiently small $\varepsilon>0$ we have $s_a(t_a + \varepsilon)<r$ and hence $a \not \sim s_a(t_a + \varepsilon)$. Because the string $s_a(t_a + \varepsilon)$ has to be associated to some row, it has to be a different row than $a$. As there are only three other rows, by taking $\varepsilon=\frac{1}{N}$ and letting $N \to \infty$ we see that some row $b$ repeats infinitely often. By taking the limit in the inequality $| \langle b, s_a\left ( t_a + \frac{1}{N} \right ) \rangle | \geq r$ we get that $| \langle b, s_a(t_a) \rangle | \geq r$ and hence the row $b$ is associated to the string $s_a(t_a)$.

Without loss of generality we can assume that $t_x = \min \{t_x, t_y, t_z, t_w \}$. We have just proved that the specific string $s_x(t_x)$ of $x$ is associated also to some other row than $x$. For the sake of simplicity we suppose that this row is $y$, but calculations are analogous for the other cases. Since  $y \sim s_x(t_x)$ we have
$$\left | \frac{1-t_x}{3}y_1 + \frac{1-t_x}{3}y_2 + \frac{1-t_x}{3}y_3 \pm t_xy_4 \right | \geq r.$$
From the triangle inequality and $y_4>0$ it follows now that
\begin{equation}
\label{nier1}
\frac{1-t_x}{3}|y_1 + y_2| + \frac{1-t_x}{3}|y_3| + t_xy_4 \geq r.
\end{equation}
Since $t_x \leq t_y$, we have $y \sim s_y(t_x)$. Thus
$$\left | \frac{1-t_x}{3}y_1 - \frac{1-t_x}{3}y_2 - t_x |y_3| + \frac{1-t_x}{3} y_4 \right | \geq r.$$
Considering the fact that $y_1, y_4>0$ and $y_2<0$, we have
$$\frac{1-t_x}{3}y_1 - \frac{1-t_x}{3}y_2 - t_x |y_3| + \frac{1-t_x}{3} y_4 \geq -t_x|y_3| \geq -\frac{1}{2} > -r$$
and hence
\begin{equation}
\label{nier2}
\frac{1-t_x}{3}y_1 - \frac{1-t_x}{3}y_2 - t_x |y_3| + \frac{1-t_x}{3} y_4 \geq r.
\end{equation}
Summation of inequalities (\ref{nier1}) and (\ref{nier2}) yields
$$\frac{1-t_x}{3}(y_1-y_2 + |y_1+y_2|) - \left ( t_x - \frac{1-t_x}{3} \right )|y_3| + \left ( \frac{1-t_x}{3}+t_x \right )y_4 \geq 2r.$$
We note that $t_x - \frac{1-t_x}{3} = \frac{4t_x - 1}{3} \geq 0$, $y_4 \leq 1$ and $y_1-y_2 + |y_1+y_2| \leq 2 \max\{|y_1|, |y_2|\} \leq 2$. Hence
$$1 = \frac{2(1-t_x)}{3}+ \left ( \frac{1-t_x}{3} + t_x \right ) \geq 2r > 1,$$
which yields the desired contradiction and finishes the first case. It is straightforward to check that when $y$ is replaced by $z$ or $w$, the proof can be carried out in the same way (albeit with some changes of coordinates and signs), so we omit the details.

We are left with the case when the matrix of $T$ can be represented as
$$\begin{bmatrix}
    + & + & + & \circ \\
    + & +  & -  & \circ \\
    + & -  & \circ  & +   \\
    + & -  & \circ  & - 
    \end{bmatrix}$$
In this case, we shall proceed in essentially the same way, but we define the functions $s_x, s_y, s_z, s_w$ a little bit differently. This time we define a function $s_x: \left [ 0, \frac{1}{2} \right ] \to \bd \mathcal{C}_4^*$ as
$$s_x(t) = \left (\frac{1}{4}, \frac{1}{4}, \frac{1}{2}-t, -\sgn(x_4) t \right ),$$
for $t \in \left [ 0, \frac{1}{2} \right ]$. From the fact that $x \sim \left ( \frac{1}{4}, \frac{1}{4}, \frac{1}{4}, \pm \frac{1}{4} \right ) \in \mathcal{S}$ it follows that $\left \langle x, s_x \left ( \frac{1}{4} \right ) \right \rangle \geq r$. On the other hand
$$0 \leq \left \langle x, s_x \left ( \frac{1}{2} \right ) \right \rangle = \frac{1}{4}(x_1 + x_2) - \frac{1}{2}|x_4| \leq \frac{1}{2}<r.$$
Therefore, there exists a unique number $t_x \in \left [ \frac{1}{4}, \frac{1}{2} \right )$ satisfying $\left \langle x, s_x \left ( t_x \right ) \right \rangle = r.$ Similarly like before, we will call the string $s_x(t_x)$ the \emph{specific string} of the row $x$. Again we have $s_x(t) \geq r \Leftrightarrow t  \leq t_x$. The functions $s_y, s_z, s_w: \left [ 0, \frac{1}{2} \right ] \to \bd \mathcal{C}_4^*$ are defined now as:
$$s_y(t) = \left (\frac{1}{4}, \frac{1}{4}, t-\frac{1}{2}, -\sgn(y_4) t \right ),$$
$$s_z(t) = \left ( \frac{1}{4}, -\frac{1}{4}, -\sgn(z_3)t, \frac{1}{2}-t \right ),$$
$$s_w(t) = \left ( \frac{1}{4}, -\frac{1}{4}, -\sgn(w_3)t, t-\frac{1}{2} \right ). $$
The numbers $t_y, t_z, t_w$ are unique numbers in $[ \frac{1}{4}, \frac{1}{2})$ such that
$$s_y(t_y)=s_z(t_z)=s_w(t_w)=r$$
and $s_y(t_y), s_z(t_z), s_w(t_w)$ are called the specific strings of the rows $y, z, w$ respectively.

Exactly like in the previous case we can prove that for every row $a \in \{x, y, z, w\}$, the specific string $s_a(t_a)$ of $a$ is associated also to some row $b \neq a$. Moreover, without loss of generality we can assume that $t_x = \min \{t_x, t_y, t_z, t_w\}$. Let $a \neq x$ be a row such that $a \sim s_x(t_x)$. Here have two essentially different possibilities to consider:
\begin{itemize}
    \item The row $a$ has its minimal element in the same column as $x$. In other words, $a=y$.
    \item The row $a$ has its minimal element in a different column that $x$. In other words, $a \in \{z, w\}$. In this case we shall assume that $a=z$, as the other case is analogous.
\end{itemize}

We start with the first case, that is $y \sim s_x(t_x)$. We have
\begin{equation}
\label{nier3}
\left | \frac{1}{4}y_1 + \frac{1}{4}y_2 + \left ( \frac{1}{2} - t_x \right )y_3  - \sgn(x_4)t_xy_4 \right | \geq r.
\end{equation}
Moreover, since $t_x \leq t_y$ we also have $y \sim s_y(t_x)$ and hence
\begin{equation}
\label{nier4}
\left | \frac{1}{4}y_1 + \frac{1}{4}y_2 + \left ( t_x - \frac{1}{2} \right )y_3  - \sgn(y_4)t_xy_4 \right | \geq r.
\end{equation}
Suppose first that $\sgn(x_4) = - \sgn(y_4)$. Then inequalities (\ref{nier3}) and (\ref{nier4}) can be restated as $|A+B|, |A-B| \geq r$, where 
$$A = \frac{y_1+y_2}{4}, \quad  B = \left ( \frac{1}{2} - t_x \right )y_3  - \sgn(x_4)t_xy_4.$$
Hence
$$2 \max \{|A|, |B|\} = |A+B| + |A-B| \geq 2r > 1.$$
However
$$|A|=\frac{y_1+y_2}{4} \leq \frac{1}{2}$$
and
$$|B| = \left | \left ( \frac{1}{2} - t_x \right )y_3  - \sgn(x_4)t_xy_4 \right | \leq \left ( \frac{1}{2} - t_x \right )|y_3| + t_x|y_4| \leq \frac{1}{2} - t_x + t_x \leq \frac{1}{2},$$
which gives a contradiction. Similarly, if $\sgn(x_4)=\sgn(y_4)$, then $|A+B|, |A-B| \geq r$, where 
$$A = \frac{y_1+y_2}{4} - \sgn(x_4)t_xy_4, \quad B = \left ( \frac{1}{2} - t_x \right )y_3.$$
Hence
$$2 \max \{|A|, |B|\} = |A+B| + |A-B| \geq 2r > 1.$$
Because all the numbers: $y_1$, $y_2$, $\sgn(x_4)t_xy_4$ are non-negative, we have
$$|A| \leq \max \left \{ \frac{y_1+y_2}{4}, t_x|y_4| \right  \} \leq \frac{1}{2}.$$
Furthermore
$$|B| = \left ( \frac{1}{2} - t_x \right )|y_3| \leq \frac{1}{2}.$$
This gives the desired contradiction and finishes the proof in the case $a=y$.

We are left with the case $a=z$, as the case $a=w$ is completely analogous. Let us suppose that $z \sim s_x(t_x)$. From the fact $t_x \leq t_z$ we know also that $z \sim s_z(t_x)$. Hence, the following inequalities are true:
\begin{equation}
\label{nier5}
\left | \frac{1}{4}z_1 + \frac{1}{4}z_2 + \left ( \frac{1}{2} - t_x \right )z_3  - \sgn(x_4)t_xz_4 \right | \geq r.
\end{equation}
\begin{equation}
\label{nier6}
\left | \frac{1}{4}z_1 - \frac{1}{4}z_2  - t_x|z_3| + \left ( \frac{1}{2} - t_x \right )z_4 \right | \geq r.
\end{equation}
First we note that the absolute value on the left-hand side of inequality (\ref{nier6}) can be omitted. Indeed, because $z_1, z_4>0$ and $z_2 < 0$ we have
$$\frac{1}{4}z_1 - \frac{1}{4}z_2  - t_x|z_3| + \left ( \frac{1}{2} - t_x \right )z_4 \geq  -t_x|z_3| \geq -\frac{1}{2} > -r.$$
Thus inequality (\ref{nier6}) rewrites as
\begin{equation}
\label{nier7}
 \frac{1}{4}(z_1-z_2)  - t_x|z_3| + \left ( \frac{1}{2} - t_x \right )z_4  \geq r.
\end{equation}
Combining inequality (\ref{nier5}) with the triangle inequality we get
\begin{equation}
\label{nier8}
\frac{1}{4}\left| z_1 + z_2 \right| + \left ( \frac{1}{2} - t_x \right )|z_3|  + t_xz_4  \geq r.
\end{equation}
Hence summation of (\ref{nier7}) and (\ref{nier8}) yields
$$\frac{1}{4} \left (z_1 - z_2 + \left| z_1 + z_2 \right| \right ) + \left ( \frac{1}{2} - 2t_x \right )|z_3| + \frac{1}{2}z_4 \geq 2r > 1.$$
However, we have also that $z_1 - z_2 + |z_1+z_2| \leq 2$, $\frac{1}{2} - 2t_x \leq 0$ and $z_4 \leq 1$. Therefore
$$\frac{1}{4} \left (z_1 - z_2 + \left| z_1 + z_2 \right| \right ) + \left ( \frac{1}{2} - 2t_x \right )|z_3| + \frac{1}{2}z_4 \leq \frac{2}{4} + 0 + \frac{1}{2} = 1.$$
We have obtained the desired contradiction and the proof is finished. \qed 

It is not clear how this four-dimensional argument could be generalized to higher dimensions. In this paper we do not focus on the general case, but in the remark below we provide an observation concerned with the asymptotic lower bound on $d_{BM}(\mathcal{C}_n, \mathcal{C}^*_n)$.

\begin{remark}
\label{remark}
Xue has conjectured that $d_{BM}(\mathcal{C}_n, \mathcal{C}^*_n) \geq \sqrt{\frac{n}{2}}$ for any $n \geq 2$ (see Conjecture 5.1 in \cite{xue}). This conjecture actually follows immediately from the well-known result of Szarek \cite{szarek}, who proved that $\frac{1}{\sqrt{2}}$ is the best possible constant in one of the variants of the Khinchin inequality. Proof of Szarek was later simplified by Haagerup \cite{haagerup}. It should be noted however, that an asymptotically better lower bound on $d_{BM}(\mathcal{C}_n, \mathcal{C}^*_n)$ is known from at least $1960$, even if not stated explicitly in the literature. It follows from some basic properties of the so-called \emph{absolute projection constant} $\lambda(X)$ of a normed space $X$. In the language of normed spaces we have $d_{BM}(\mathcal{C}_n, \mathcal{C}^*_n) = d_{BM}(\ell_{\infty}^n, \ell_{1}^n)$. It is widely known that $\lambda(\ell_{\infty}^n)=1$ and already in 1960 Gr\"unbaum \cite{grunbaum2} has determined the absolute projection constant of the space $\ell_{1}^n$. From his result it follows that $\frac{\lambda(\ell_1^n)}{\sqrt{n}} \to \sqrt{\frac{2}{\pi}}$ as $n \to \infty$. Combining this with the well-known inequality $d_{BM}(X, Y) \geq \frac{\lambda(X)}{\lambda(Y)}$ (true for any $n$-dimensional normed spaces $X, Y$, see Corollary 6 in Section III.B. in \cite{wojtaszczyk}), we obtain an asymptotic lower bound  $d_{BM}(\mathcal{C}_n, \mathcal{C}^*_n) \geq \lambda(\ell_1^n) \sim\sqrt{\frac{2}{\pi}n}.$ From the viewpoint of asymptotics, this seems to be the best lower bound  currently known. It is not clear however, if the constant $\sqrt{\frac{2}{\pi}}$ is asymptotically the best possible.

\end{remark}

\section{Planar convex bodies equidistant to symmetric convex bodies}
\label{sect2d}
In this section we establish a large family of planar convex bodies that are equidistant to the whole family of symmetric convex bodies. It is well known that the triangle is equidistant to all symmetric convex bodies with the distance equal to $2$. Our construction shows that there are much more planar convex bodies with this property than just a triangle. In particular, for each $r \in \left (\frac{7}{4}, 2 \right )$ there are continuum many affinely non-equivalent convex pentagons equidistant to symmetric convex bodies with the distance $r$. Our main tool is a classical concept of the convex geometry: the \emph{asymmetry constant}. For a given convex body $K \subseteq \mathbb{R}^2$ we define its asymmetry constant $\as(K)$ as
$$\as(K) = \inf \{ r>0 : \text{ there exists } z \in \inte K \text{ such that } K-z \subseteq -r(K-z) \}.$$
In the planar case it is known that there exists exactly one point $z \in \inte K$ for which this infimum is attained. Such a point $z$ is called a \emph{Minkowski center} of $K$. The following properties of the asymmetry constant are well-known for a convex body $K \subseteq \mathbb{R}^2$ (see for example \cite{grunbaum}, \cite{neumann}):

\begin{enumerate}
\item $1 \leq \as(K) \leq 2$,
\item $\as(K)=1$ if and only if $K$ is symmetric,
\item $\as(K)=2$ if and only if $K$ is a triangle,
\item If $z$ is the Minkowski center of $K$, than the boundaries of convex bodies: $K-z$ and $-\as(K)(K-z)$ intersect in at least three points.
\end{enumerate}

The asymmetry constant relates to the Banach-Mazur distance in the following natural way. The result is folklore, for a short proof see for example Proposition 3.1 in \cite{brandeberg}.

\begin{lem}
\label{asbm}
For every convex body $K \subseteq \mathbb{R}^n$ we have
$$\as(K) = \inf d_{BM}(K, L),$$
where the infimum runs over all symmetric convex bodies $L \subseteq \mathbb{R}^n$.
\end{lem}

From the properties above it follows immediately that if $S \subseteq \mathbb{R}^2$ is a triangle, then $d_{BM}(S, L) \geq 2$ for any symmetric convex body $L$. The opposite inequality follows from a classical maximal area argument -- it is easy to prove that if $S \subseteq L$ is a triangle with the maximal possible area, then $L$ is contained in a copy of $S$ scaled by $2$. In our construction we will proceed in a very similar way. The lower bound will follow from the asymmetry constant and Lemma \ref{asbm}, while for the upper bound we will use a triangle of the maximal area.

By $u_1, u_2, u_3$ we denote the vertices of an equilateral triangle in $\mathbb{R}^2$. For a standard scalar product $\langle \cdot, \cdot \rangle$ in $\mathbb{R}^2$ we assume that $\langle u_i, u_i \rangle = 1$ for each $1 \leq i \leq 3$ and $\langle u_i, u_j \rangle = -\frac{1}{2}$ for $i \neq j$. In particular we have $u_1+u_2+u_3=0$.

\begin{figure}
    \centering
    \includegraphics[scale=0.2]{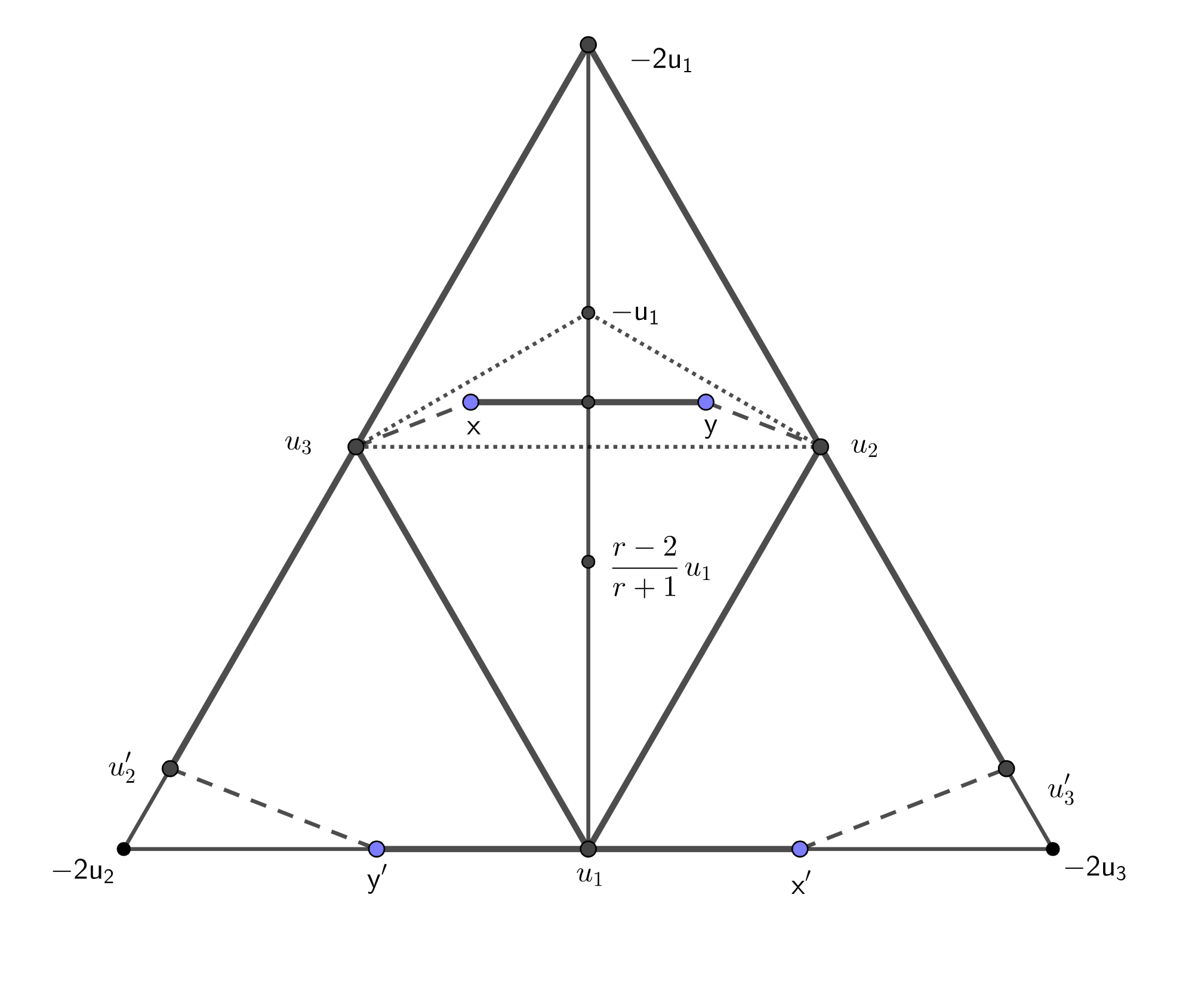}
    \caption{Construction of an axially symmetric convex body with the distance $r$ to every symmetric convex body (in the picture $r=1.8$). A convex pentagon $\conv\{u_1, u_2, u_3, x, y\}$ satisfies this condition. If $r$ is in the open interval $\left ( \frac{7}{4}, 2 \right )$ and $k$ is in the interval $\left [\frac{1}{2r}, 2 - \frac{3}{r} \right )$, then $x$ and $u_3$ instead of being connected by the dashed segment, could be joined by any convex curve that is inside the triangle $\conv\{u_2, u_3, -u_1\}$. The point $\frac{r-2}{r+1}u_1$ is the Minkowski center of $K$. Points $u_2', u_3', x', y'$ are corresponding points in a homothetical image of $K$ with ratio $-r$, while the corresponding point for $u_1$ is $-2u_1$.}
    \label{ryskonstr}
\end{figure}

The following lemma contains our main construction. It should be noted that the inequality $k \leq 2- \frac{3}{r}$ in the second condition guarantees that the first two conditions are not excluding each other. In fact, for $k = 2 - \frac{3}{r}$ the endpoints of the given segment belong to the sides of the quadrilateral $\conv \{-u_1, u_1, u_2, u_3\}$ and for smaller $k$ they lie in its interior. See Figure \ref{ryskonstr} for an illustration. 

\begin{lem}
\label{lemkonstr}
Let $\frac{7}{4} \leq r \leq 2$. Suppose that a convex body $K \subseteq \mathbb{R}^2$ satisfies the following conditions:
\begin{enumerate}
\item $\conv\{u_1, u_2, u_3\} \subseteq K \subseteq \conv \{-u_1, u_1, u_2, u_3\}$.
\item The boundary of $K$ contains a segment
$$\left [ \left ( \frac{r-3}{r} + k \right )u_1 + 2ku_2,  \left ( \frac{r-3}{r} + k \right )u_1 + 2ku_3\right ],$$
where $k$ is a fixed real number in the interval  $\left [\frac{1}{2r}, 2 - \frac{3}{r} \right ]$.
\item The line $\{ x \in \mathbb{R}^2 : \langle x, u_2 \rangle = \langle x, u_3 \rangle \}$ is a symmetry axis of $K$.
\end{enumerate}
Then 
$$d_{BM}(K, L)=\as(K)=r$$
for every symmetric convex body $L \subseteq \mathbb{R}^2$.
\end{lem}
\emph{Proof.} We denote $S = \conv\{u_1, u_2, u_3\}$. Since for planar convex bodies the Minkowski center is unique and $K$ has a symmetry axis, the Minkowski center of $K$ is a point of the form $\alpha u_1$, where $\alpha \in \mathbb{R}$. To determine $\alpha$ we note that by assumption, the line passing through $-2u_1, -2u_2$ and the line passing through $u_1, u_2$ are two different parallel lines supporting $K$. Hence, the homothety with center $\alpha u_1$ and ratio $-\frac{1}{as(K)}$ sends $u_3$ (lying on the first line) to some point lying on the line through $u_1, u_2$, which can be described as  $\{ x: \ \langle x, u_3 \rangle = -\frac{1}{2} \}$. The image of $u_3$ in this homothety is equal to
$$\frac{1}{\as(K)}\left ((1+\as(K))\alpha u_1 - u_3 \right )$$
and hence
$$-\frac{1}{2}=\frac{1}{\as(K)} \left \langle\left ((1+\as(K))\alpha u_1 - u_3 \right ), u_3 \right \rangle = -\frac{1}{\as(K)} \left ( \frac{(1+\as(K))\alpha}{2} +1 \right ),$$
which yields the equality
$$\alpha = \frac{\as(K)-2}{\as(K)+1}.$$
On the other hand, by the assumption $K$ also has two lines parallel to $u_2u_3$ in the boundary. In consequence the homothety with center $\alpha u_1$ and ratio $-\frac{1}{\as(K)}$ sends $u_1$ to $\frac{r-3}{r}u_1$. By a direct calculation we get the following
\begin{equation}
\label{minkowski}
\as(K)=r \text{ and the point } \frac{r-2}{r+1}u_1 \text{ is the Minkowski center of K}.
\end{equation}
Let us denote by $K'$ the homothetical image of $K$ with center $\frac{r-2}{r+1}u_1$ and ratio $-r$. It is now also easy to verify that in this homothety the image of the point $u_1$ is equal to $-2u_1$. Thus the convex body $K'$ contains a parallelogram $\conv \{-2u_1, u_1, u_2, u_3\}$.


Now, let $L \subseteq \mathbb{R}^2$ be any symmetric convex body. Our goal is to find an affine image $L_0$ of $L$ such that 
$$K \subseteq L_0 \subseteq K'.$$
Indeed, if an affine copy $L_0$ of $L$ satisfies this inclusions, then for a certain $u$ we have $L_0 \subseteq -rK + u$ or $-L_0 \subseteq rK - u$. If $s$ is the center of symmetry $L_0$, then $L_0 = 2s-L_0$ and hence $L_0 \subseteq rK + (2s-u)$. Thus $L_0$ is contained between two homothetical copies of $K$ with the ratio $r$.

In order to prove the inclusions above, let us consider a triangle $\conv \{a, b, c\} \subseteq L$ with a maximal possible area among all triangles contained in $L$. Clearly, the symmetry center of $L$ lies in the triangle $abc$, as otherwise we could easily find a triangle contained in $L$ with a larger area. Let $g$ be the center of gravity of the triangle $abc$. The triangle $abc$ is divided into three triangles: $gab, gbc$ and $gca$. Without loss of generality, we assume that the center of symmetry of $L$ lies in the triangle $\{g, b, c\}$. Now we consider an affine transformation $T: \mathbb{R}^2 \to \mathbb{R}^2$ defined by conditions: $T(a)=u_1$, $T(b)=u_2$, $T(c)=u_3$ and we denote $L_0=T(L)$. We shall prove that $L_0$ is the desired affine image of $L$.

From the fact that the triangle $S$ is of the maximal area in $L_0$, it follows that the line passing through $u_1$ and parallel to $u_2u_3$ is supporting $L_0$. Similarly for $u_2$ and $u_3$. Hence we have that $L_0 \subseteq -2S$. We start with proving the inclusion $K \subseteq L_0$.

\begin{figure}
    \centering
    \includegraphics[scale=0.15]{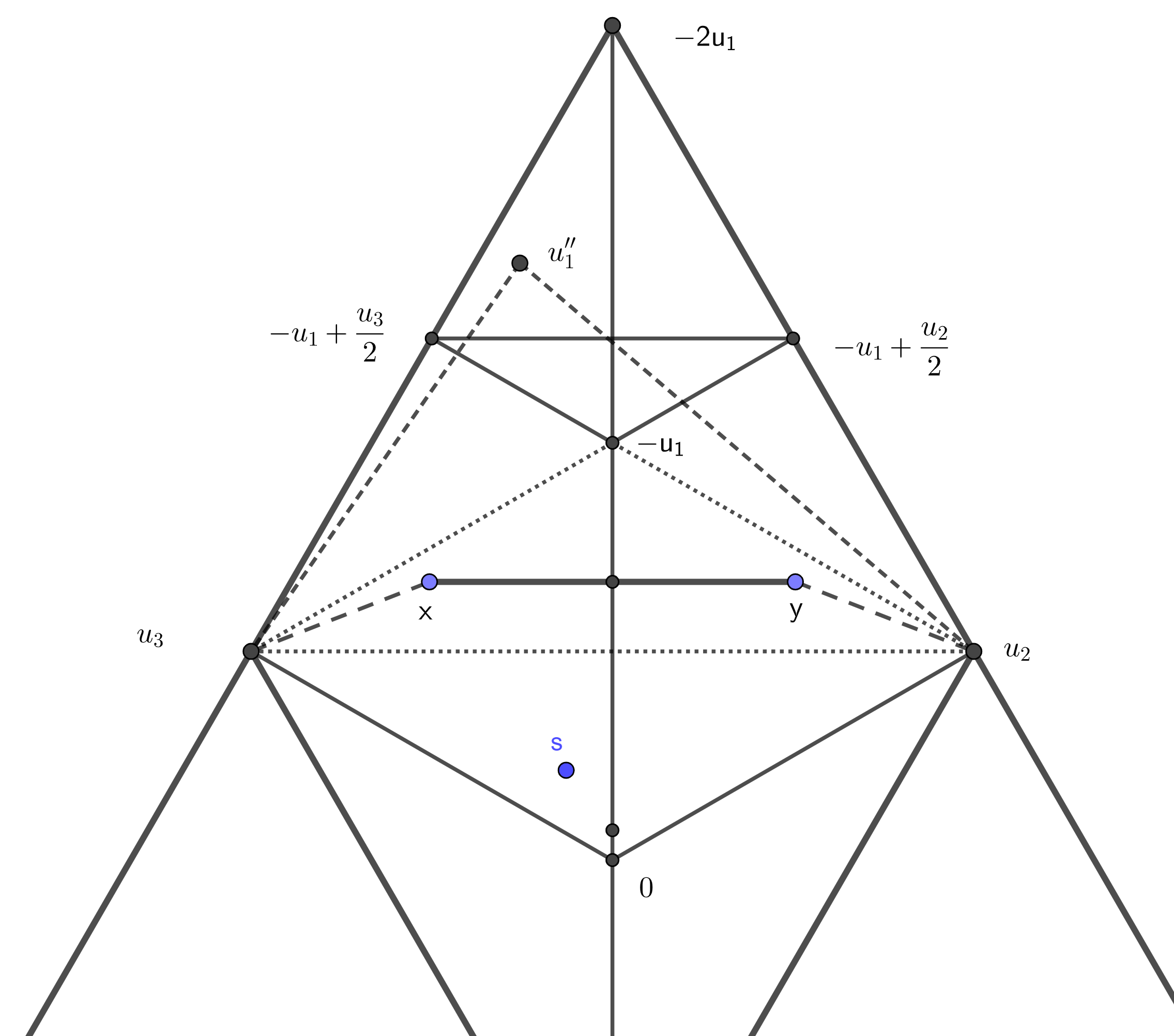}
    \caption{Proof of the inclusion $K \subseteq L_0$. By the assumption the symmetry center $s$ of $L_0$ is inside the triangle with vertices $0, u_2, u_3$. The reflection $u_1''=2s-u_1$ of $u_1$ lies in the quadrilateral with vertices $-u_1, -u_1 + \frac{u_2}{2}, -u_1 + \frac{u_3}{2}, -2u_1$. It follows that the point $-u_1$ is inside the triangle $u_1''u_2u_3$.}
    \label{rysinclusion}
\end{figure}

By the assumption we have that  $K \subseteq \conv \{-u_1, u_1, u_2, u_3\}$. Thus it is enough to check that $-u_1 \in L_0$. Let $s \in \conv \{0, u_2, u_3\}$ be the symmetry center of $L_0$. The reflection $u_1''=2s-u_1$ of $u_1$ lies in $L_0$. Because $s$ lies in the triangle $\conv\{0, u_2, u_3\}$, the reflection $u_1''$ belongs to the triangle $S'$ with vertices $-u_1, 2u_2-u_1, 2u_3-u_1$. However, because $L_0 \subseteq -2S$, point $u_1''$ belongs to the intersection $(-2S) \cap S'$, which is a quadrilateral with vertices $-u_1, -u_1 + \frac{u_2}{2}, -u_1 + \frac{u_3}{2}, -2u_1$. It is now clear, that regardless of the position of $u_1''$ inside this quadrilateral, the triangle with the vertices $u_1'', u_2, u_3$ contains $-u_1$ (see Figure \ref{rysinclusion}) and it follows from the convexity of $L_0$ that $-u_1 \in L_0$. This concludes the proof of the first inclusion.

Now we shall prove that $L_0 \subseteq K'$. Because the convex body $K'$ contains the whole parallelogram $\conv\{u_1, u_2, u_3, -2u_1\}$, it is enough to check that $L_0 \cap \conv\{u_1, u_3, -2u_2\} \subseteq K' $ and $L_0 \cap \conv\{u_1, u_2, -2u_3\} \subseteq K' $. We will check the first inclusion, as the second one can be verified in a completely analogous way.

Let $u_2'' = 2s - u \in L_0$ be the reflection of $u_2$. To show that $L_0 \cap \conv\{u_1, u_2, -2u_3\} \subseteq K'$ we will establish the following inequality
\begin{equation}
\label{ilskal}
\langle u_2'', u_2 \rangle \geq -\frac{5}{4}
\end{equation}
We have assumed that $s \in \{0, u_2, u_3\}$, so let us write $s = Au_2 + Bu_3$, where $A, B \geq 0$ and $A+B \leq 1$. Then
$$u_2'' = 2s - u_2 = (2A-1)u_2 + 2Bu_3.$$
Since $u_2'' \in L_0$ and $L_0 \subseteq -2S$, we have that $u_2'' \in -2S$. We can assume that $u_2''$ belongs to the triangle with vertices $u_1, u_3, -2u_2$, as for every point $u \in (-2S) \setminus \conv\{u_1, u_3, -2u_2\}$ we have that $\langle u, u_2 \rangle \geq -\frac{1}{2}$ (and thus inequality (\ref{ilskal}) is satisfied). Hence let us write
$$u_2'' = (-2E)u_2 + Fu_3 + (1-E-F)u_1 = (-2E)u_2 + Fu_3 + (E+F-1)(u_2+u_3)$$
$$ =(F-E-1)u_2 + (E+2F-1)u_3,$$
where $E, F \geq 0$ and $E+F \leq 1$. Therefore we have that
$$(2A-1)u_2 + 2Bu_3 = (F-E-1)u_2 + (E+2F-1)u_3$$
and consequently
$$2A = F - E \quad \text{ and } \quad 2B = E+2F-1.$$
In particular $F \geq E$ and hence $E \leq \frac{1}{2}$. Thus
$$\langle u_2'', u_2 \rangle = \left \langle (F-E-1)u_2 + (E+2F-1)u_3, u_2 \right  \rangle = F-E-1 + - \frac{E}{2} - F + \frac{1}{2}$$
$$ =-\frac{1}{2}\left (1+3E \right ) \geq -\frac{1}{2} \left ( 1 + \frac{3}{2} \right ) = -\frac{5}{4},$$
which proves inequality (\ref{ilskal}).

\begin{figure}
    \centering
    \includegraphics[scale=0.15]{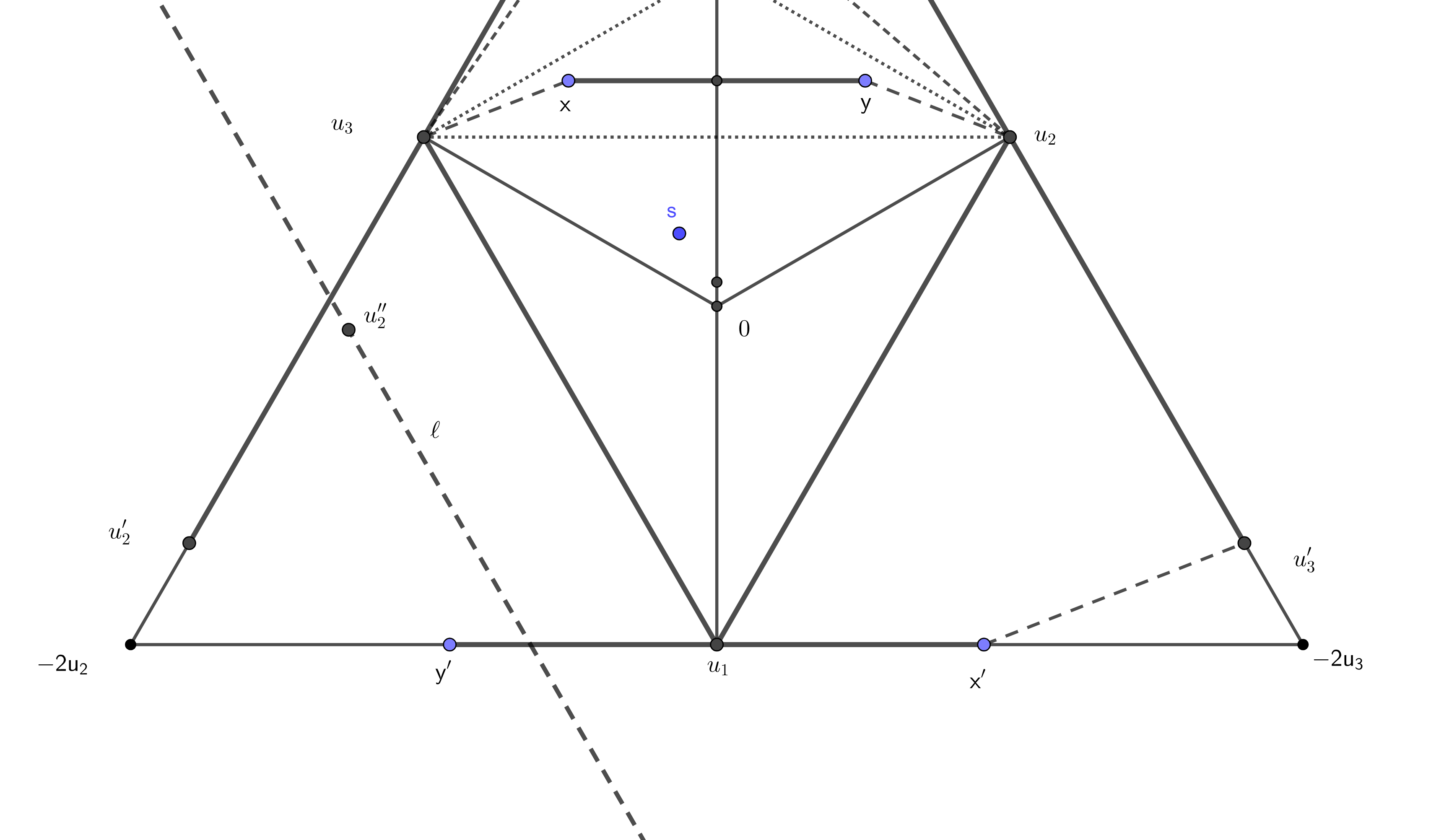}
    \caption{Proof of the inclusion $L_0 \cap \conv \{u_1, u_3, -2u_2\} \subseteq K'$. Because the line $\ell$ is supporting to the convex body $L_0$ at the point $u_2''$, the convex body $L_0$ is on the right side of $\ell$ in the picture. The points $u_2'$ and $y'$ are to the left of $\ell$.}
    \label{rysinclusion2}
\end{figure}

To establish the inclusion $L_0 \cap \conv \{u_1, u_3, -2u_2\} \subseteq K'$, let us consider a supporting line $\ell$ to $L_0$ at $u_2''$, parallel to the line $u_1u_3$. Line $\ell$ is of the form $\{x \in \mathbb{R}^2: \ \langle x, u_2 \rangle = \gamma\}$, where $\gamma \geq -\frac{5}{4}$ by inequality (\ref{ilskal}). To prove the desired inclusion we shall show that $K'$ has two points on two sides of the triangle $-2S$, which are further to the ,,left'' of the line $\ell$ in the direction of $u_2$ (see Figure \ref{rysinclusion2}). More formally, it is enough to check that the scalar product with $u_2$ of these two points of $K'$ is not greater than $-\frac{5}{4}.$

The point $y = \left ( \frac{r-3}{r} + k \right )u_1 + 2ku_2$ is in $K$ by the assumption. The corresponding point $y'$ of $y$ in $K'$ lies on the side $[-2u_2, -2u_3]$ of $-2S$. To calculate it explicitly, we use the property (\ref{minkowski}):
$$y' = (r+1)\left (\frac{r-2}{r+1}u_1 \right ) - ry = (r-2)u_1 - (r-3 - rk)u_1 - 2rk u_2$$
$$=(1-rk)u_1 - 2rk u_2.$$
Hence
$$\langle y', u_2 \rangle = \frac{rk-1}{2} - 2rk = -\frac{1}{2} (3rk + 1) \leq - \frac{1}{2} \left (\frac{3}{2} + 1 \right ) = - \frac{5}{4}.$$
Similarly, the corresponding point $u_2'$ of $u_2 \in K$ in $K'$ is given by
$$u_2'=(r+1)\left (\frac{r-2}{r+1}u_1 \right ) - ru_2 = (r-2)u_1 - ru_2.$$
This point lies on the side $[-2u_2, -2u_1]$ of $-2S$ and
$$\langle u_2', u_2 \rangle=\langle (r-2)u_1 - ru_2, u_2 \rangle = \frac{2-r}{2} - r = 1 - \frac{3}{2}r \leq 1 - \left (\frac{3}{2} \cdot \frac{7}{4} \right ) = - \frac{13}{8} <  - \frac{5}{4}.$$
In this way we have verified that $K'$ has two points on two sides of $-2S$ with the scalar product with $u_2$ not greater than $-\frac{5}{4}$. From convexity of $K'$ it follows that $L_0 \cap \{u_1, u_3, -2u_2\} \subseteq K'$ and the proof is finished. \qed 

In the next theorem we summarize our results about convex bodies equidistant to the symmetric bodies.

\begin{twr}
\label{twr2d}
For every $\frac{7}{4} \leq r \leq 2$ there exists a convex pentagon $K \subseteq \mathbb{R}^2$ satisfying $d_{BM}(K, L)=\as(K)=r$ for every symmetric convex body $L \subseteq \mathbb{R}^2$. If $r>\frac{7}{4}$, then there are continuum many affinely non-equivalent convex pentagons $K$ with this property.

Moreover, if a convex body $K \subseteq \mathbb{R}^2$ has this property for some $r$, then $r \geq \sqrt{\frac{3}{2}}$ and $K$ is not smooth and not strictly convex.
\end{twr}

\emph{Proof.} Let us start with the first part. Directly from Lemma \ref{lemkonstr} it follows that for $\frac{7}{4} \leq r < 2$ a convex pentagon 
$$K = \conv \left \{u_1, u_2, u_3, \left ( \frac{r-3}{r} + k \right )u_1 + 2ku_2,  \left ( \frac{r-3}{r} + k \right )u_1 + 2ku_3 \right \},$$
where $k = 2 - \frac{3}{r}$, satisfies the desired conditions. Moreover, if $r>\frac{7}{4}$, then we have continuum many possibilities for $k \in \left ( \frac{1}{2r}, 2 - \frac{3}{r} \right ) $. It is easy to see that if $k_1 \neq k_2$, then the convex pentagons corresponding to $k_1$ and $k_2$ are not affinely equivalent. Indeed, both of them has exactly one pair of side and diagonal that are parallel to each other. It follows that these pairs have to be mapped to each other by any affine transformation mapping one pentagon to the other. However, the ratio of lengths of parallel segments remains the same in any affine mapping, but these ratios are clearly different if $k_1 \neq k_2$. 

For the second part we note, that it is known that the Banach-Mazur distance between the square and the regular hexagon is equal to $\frac{3}{2}$ (see \cite{lassak2}). Therefore, if a convex body $K$ is of distance $r$ to both of them, then by triangle inequality we clearly have $r^2 \geq \frac{3}{2}$. To establish the second part, we will prove actually a much more general fact: if a convex body $K \subseteq \mathbb{R}^2$ satisfies $d_{BM}(K, \mathcal{C}_2)=\as(K)$, then $K$ is not smooth and not strictly convex.

Let us suppose that a convex body $K \subseteq \mathbb{R}^2$ satisfies $d_{BM}(K, \mathcal{C}_2)=r$, where $r = \as(K)$. Without loss of generality we can assume that $0$ is the Minkowski center of $K$. If $r=1$, then $K$ is a parallelogram and there is nothing to prove. Hence we can assume that $K$ is not centrally symmetric. We will rely on the following well-known fact: the boundaries of the convex bodies $K$ and $-rK$ have at least three points of contact that are not all on one line (thus forming a triangle $\mathcal{T}_1$) and there exist common supporting lines to $K$ and $-rK$ at these three points that form a triangle $\mathcal{T}_2$ containing both $K$ and $-rK$. This was established in a classical paper of B.H Neumann \cite{neumann} (see Sections 3 and 4). Therefore we have a chain of inclusions $\mathcal{T}_1 \subseteq K \subseteq -rK \subseteq \mathcal{T}_2$.

If $d_{BM}(K, \mathcal{C}_2)=r$, then there exists a parallelogram $\mathcal{P} \subseteq \mathbb{R}^2$ such that 
$$K \subseteq \mathcal{P} \subseteq rK + v,$$
for some vector $v \in \mathbb{R}^2$. However, if $s$ is the center of symmetry of $\mathcal{P}$, then $2s - \mathcal{P} = \mathcal{P}$ and hence
$$-\mathcal{P} \subseteq rK + (v-2s),$$
which yields
$$\mathcal{P} \subseteq -rK + (2s-v).$$
Thus we have $K \subseteq -rK + (2s-v)$. However, because in the planar case the Minkowski center of $K$ is unique, we must have $2s=v$. It follows that $\mathcal{P} \subseteq -rK$. To summarize, we have the following chain of inclusions
$$\mathcal{T}_1 \subseteq K \subseteq \mathcal{P} \subseteq -rK \subseteq \mathcal{T}_2.$$
Vertices of the smaller triangle $\mathcal{T}_1$ are on the sides of the large triangle $\mathcal{T}_2$ and hence they lie on the boundaries of $K$, $-rK$ and $\mathcal{P}$. We start with observing that no two vertices of $\mathcal{T}_1$ can lie in the interiors of some opposite sides of $\mathcal{P}$. Indeed, for any point lying in the interior of a side of $\mathcal{P}$, there exists a unique supporting line to $\mathcal{P}$ at this point -- namely the line determined by this side. However, from the inclusions $\mathcal{T}_1 \subseteq \mathcal{P} \subseteq \mathcal{T}_2$ it follows that lines determined by the sides of $\mathcal{T}_2$ are supporting at the vertices of $\mathcal{T}_1$ to $\mathcal{P}$ and no two of them are parallel.

Thus at least one vertex $x$ of $\mathcal{T}_1$ is also a vertex of $\mathcal{P}$. However, $x$ is also a boundary point of $K$ and from the inclusion $K \subseteq \mathcal{P}$ it follows that any supporting line of $\mathcal{P}$ to $x$ is also supporting line of $K$. Two lines determined by the sides of $\mathcal{P}$ containing $x$ are two different supporting lines at $x$. This shows that $K$ has at least two different supporting line at $x$ and hence $K$ is not smooth.

Now we shall prove that $K$ is not strictly convex. If there exists a vertex $x$ of $\mathcal{T}_1$ lying on the side of $\mathcal{P}$, then from the fact that $x$ is a boundary point of $\mathcal{T}_2$ and the inclusion $\mathcal{P} \subseteq \mathcal{T}_2$ it follows that this side of $\mathcal{P}$ is contained in a side of $\mathcal{T}_2$. Since $\mathcal{P} \subseteq -rK \subseteq \mathcal{T}_2$ this side is contained in the boundary of $-rK$ and we conclude that the convex body $-rK$ contains a segment in its boundary. Thus $-rK$ is not strictly convex and the same holds obviously also for $K$.

We are left with the situation, in which every vertex of $\mathcal{T}_1$ is also a vertex of $\mathcal{P}$. In this case, the convex body $K$ contains a consecutive pair of vertices of $\mathcal{P}$ in its boundary. Therefore, since $K \subseteq \mathcal{P}$ it contains the whole side of $\mathcal{P}$ in the boundary. Again we conclude that $K$ contains a segment in its boundary and thus it is not strictly convex. This finishes the proof. \qed

It should be noted that the construction given in Lemma \ref{lemkonstr} yields much more convex bodies equidistant to the symmetric bodies than just convex pentagons. For every $r \in \left ( \frac{7}{4}, 2 \right )$ and $k \in \left ( \frac{1}{2r}, 2 - \frac{3}{r} \right )$ it is possible to connect points $u_3$ and $x$ (using the notation of Figure \ref{ryskonstr}) with any convex curve lying in the triangle $\conv{u_2, u_3, -u_1}$ (the points $u_2$ and $y$ are then connected with the symmetric curve). Therefore such a convex body does not necessarily need to be a polygon, but as we have already seen, it can not be strictly smooth or convex. We do not know if convex bodies in $\mathbb{R}^n$ that are equidistant to all symmetric convex bodies and are different from a simplex exist for all $n \geq 2$, but it is highly possible. In the planar case it would be interesting to determine the smallest possible $r$, for which there exists a planar convex body with the distance $r$ to every symmetric convex body.

\end{document}